\documentclass[11pt]{article}
\usepackage{}
\usepackage{amsmath,amsfonts}
\usepackage{txfonts,wasysym,lscape,amssymb,mathrsfs,extarrows}
\usepackage[all,pdf]{xy}
\usepackage{graphicx}
\usepackage{subfigure}
 \usepackage{multirow}
 \usepackage{rotating}
 \usepackage[table,xcdraw]{xcolor}
\RequirePackage{color,graphicx}

\newtheorem{theorem}{Theorem}[section]
\newtheorem{proposition}{Proposition}[section]
\newtheorem{lemma}{Lemma}[section]

\newtheorem{algorithm}{Algorithm}[section]
\newtheorem{definition}{Definition}[section]

\newtheorem{remark}{Remark}[section]

\newtheorem{assumption}{Assumption}[section]

\title{Optimality Conditions and Numerical Algorithms for A Class of Linearly Constrained Minimax Optimization Problems}

\author{Yu-Hong Dai\footnote{LSEC, AMSS, Chinese Academy of Sciences, Beijing 100190, P. R. China. {\sl Email}: dyh@lsec.cc.ac.cn.
This author was supported by the Natural Science Foundation of China (Nos. 11991021, 11991020, 12021001 and 11971372) and
the Strategic Priority Research Program of Chinese Academy of Sciences (no. XDA27000000).}
 \footnote{School of Mathematical Sciences, University of Chinese Academy of Sciences, Beijing 100049, P. R. China.}\quad Jiani Wang\footnote{LSEC, AMSS, Chinese Academy of Sciences, Beijing 100190, P. R. China. {\sl Email}: wjiani@lsec.cc.ac.cn.} \quad and \quad Liwei Zhang \footnote{ School of Mathematical Sciences, Dalian University of Technology, Dalian 116024, P. R. China. {\sl Email}: lwzhang@dlut.edu.cn. This author was supported by the Natural Science Foundation of China (Nos. 11971089 and 11731013).}}

\date{}
\begin{document}

\maketitle

\begin{abstract}
It is well known that there have been many numerical algorithms for solving nonsmooth minimax problems, numerical algorithms for nonsmooth minimax problems with joint linear constraints are very rare. This paper aims to discuss optimality conditions and develop  practical numerical algorithms for minimax problems with joint linear constraints. First of all, we use the properties
of  proximal mapping and KKT system to establish optimality conditions. Secondly, we propose a framework of alternating coordinate algorithm for the minimax problem and analyze its convergence properties. Thirdly, we develop a proximal gradient multi-step ascent decent method (PGmsAD) as a numerical algorithm and demonstrate that the method can find an $\epsilon$-stationary point for this kind of nonsmooth nonconvex-nonconcave problem in ${\cal O}\left(\epsilon^{-2}\log  \epsilon^{-1}  \right)$ iterations. Finally, we apply PGmsAD to generalized absolute value equations, generalized linear projection equations and linear regression
problems and report the efficiency of PGmsAD on  large-scale optimization.
\vskip 6 true pt \noindent \textbf{Key words}: minimax optimization, proximal mapping, proximal gradient multi-step ascent decent method, iteration complexity, generalized absolute value equations, linear regression, generalized linear projection equations.
\vskip 12 true pt \noindent \textbf{AMS subject classification}: 90C30
\end{abstract}
\bigskip\noindent

\section{Problem setting}
 \setcounter{equation}{0}
The main purpose of  this paper is to build optimality conditions and numerical algorithms for the minimax optimization problem with joint
linear constraints,
\begin{equation}\label{minimax}
\begin{array}{rl}
\displaystyle \min_{x \in {\Re^n}}\max_{y \in {\Re^m}}& f(x,y)=\varphi (x)+g(x)+x^TKy-h(y)-\psi (y)\\[10pt]
\mbox{subject to} & Ax+By+c=0,
\end{array}
\end{equation}
where  $g: \Re^n \rightarrow \Re$ and $h:\Re^m \rightarrow \Re$ are real-valued smooth functions, $\varphi:\Re^n \rightarrow \overline \Re$ and $\psi:\Re^m \rightarrow \overline \Re$ are extended real-valued proper lower semicontinuous convex  functions (here $\overline \Re=\Re \cup \{+\infty\}$), $K\in\Re^{n\times m}$, $A \in \Re^{q \times n}$, $B \in \Re^{q\times m}$ and $c \in \Re^q$.

\subsection{Applications}
Minimax problems \eqref{minimax} play a significant role in optimization since many interesting linear equations can be converted to equivalent linearly constrained minimax problems. Traditional methods for solving linear equations include least squares method \cite{Ben90}, fixed point method \cite{Agarwal2001}, Newton method \cite{Akram2015}. However, these classical numerical algorithms cannot handle equations with absolute values or projection operators, such as generalized absolute value equations (GAVE) \cite{Zhou2021} and generalized linear projection equations (GLPE) \cite{Eshaghnezhad2016}. Recently, numerical algorithms have been proposed only to solve such equations with special structures, such as generalized Newton methods \cite{Mangasarian2008}, neural networks \cite{XiaWang2000}, {\it etc.} It is noteworthy that these intractable linear system problems can be translated into equivalent linearly constrained minimax problems.

\textbf{Generalized absolute value equations (GAVE)}. GAVE is a popular nonsmooth NP-hard problem in the form of
\begin{equation}\label{GAVE}
  Ax+B|x|=b,
\end{equation}
where $A,\ B\in \mathfrak{R}^{m\times n},\ b\in \mathfrak{R}^{m}$ and for $x=(x_1,\ldots,x_n)\in\mathfrak{R}^{n}$, $|x|:=(|x_1|,\ldots,|x_n|)\in \mathfrak{R}^{n}_+$. As an important tool in the field of optimization, GAVE is widely used to solve problems in diverse fields, including nonnegative constrained least squares problems, quadratic programming, complementarity problem, bimatrix games ({\it e.g.} \cite{Cottle2009,Enkhbat2021,Mangasarian2006,Zheng2016}). In this paper, we prove that GAVE is equivalent to a smooth nonspareable linearly constraint convex-concave minimax problem as follows:
\begin{equation*}
  \begin{array}{ll}
\displaystyle \min_{x\in {\Re^n_+}}\max_{z\in {\Re^m_+},y\in \Re^m}& (b-(A+B)x)^Ty\\[2pt]
\mbox{subject to} & x-(B-A)^Ty-z=0.
\end{array}
\end{equation*}
Moreover, we also focus on the well-known linear regression problems with joint linearly constraints and strongly-convex-strongly-concave quadratic objective functions.

\textbf{Generalized Linear Projection Equations (GLPE)}. As one of the important techniques for solving constrained optimization problems, GLPE has attracted extensive attention in linear variational inequalities, fixed point problems, bimatrix equilibrium points and traffic network modeling (e.g. \cite{Cottle2009,Facchinei2003,XiaWang2004}). In this subsection, we are concerned with solving the generalized linear projection
equations of the following form
\begin{equation}\label{glpe}
  Ax+Bx_K=b,\ {\rm subject\ to}\ x_K=P_K(x),
\end{equation}
where $A,\ B\in \mathfrak{R}^{m\times n},\ b\in \mathfrak{R}^{m}$, $K\subset\mathfrak{R}^{n}$ is a closed convex set and for any $x\in \mathfrak{R}^{n}$, $P_K(x)$ represents the projection of $x$ onto $K$. If $K$ is a cone,
GLPE can be translated into a nonsmooth nonspareable linearly constraint convex-concave minimax problem
\begin{equation*}
  \begin{array}{ll}
\displaystyle \min_{x\in {\Re^n}}\max_{y\in \Re^m,z\in {\Re^n}}& \delta_K(x)+(b-(A+B)x)^Ty-\delta_{K^{\circ}}(z)\\[2pt]
\mbox{subject to} & x-A^Ty-z=0,
\end{array}
\end{equation*}
where $K^{\circ}$ represents the polar cone of $K$ and $\delta_K(x)$ is an indicator function.

Motivated by the relationship between linear equations and minimax problems, we will propose numerical algorithms for linearly constrained minimax problems \eqref{minimax} and apply these algorithms to solve linear equations with absolute values or projection operators.

\subsection{Related works}
There are many effective algorithms that can be solved when $A=0$ or $B=0$ in linearly constrained minimax problems \eqref{minimax}, such as the
alternating coordinate method \cite{Arrow58}, first-order primal-dual algorithm \cite{CPock2011} and the multi-step gradient descent ascent method \cite{Raza2019}. However, when $A\neq0$ and $B\neq0$, since the variables $x,\ y$ in the constraint are coupled, \eqref{minimax} can not be solved by the above conventional algorithms. It is noted that, by introducing Lagrange multiplier and establishing Lagrange function, we can solve nonlinear programming problems effectively and quickly (see \cite{B82}). Similarly, we hope to design an effective algorithm to solve the joint linearly constrained minimax problem by establishing the Lagrange function of \eqref{minimax}. Since $\varphi$ and $\psi$ in \eqref{minimax} are nonsmooth, the proximal mapping technique is used to deal with nonsmooth terms.

The study of algorithms for solving minimax problems is always active. For the case when variables $x$ and $y$ have separable constraints, there are many publications about constructing and analyzing numerical algorithms for the minimax problem, such as \cite{Arrow58,CPock2011,CPock2016a,CPock2016b,ChenLan2014,LinJordan2020,Nem2004,Nes2007,Raza2019,Raza2018,XuLan2020}. Just recently, proximal projection gradient-type algorithms have been deeply studied, for nonsmooth minimax optimization problems. Valkonen \cite{Valkonen2014} proposed  an
extension of the modified primal-dual hybrid gradient method, due to Chambolle and Pock \cite{CPock2011}. Mokhtari {\it et al.} \cite{Ozd2019a} proposed algorithms admitting a unified analysis as approximations of the classical proximal point method for solving saddle point problems. Dai {\it et al.} \cite{DaiWZ2020} proposed a semi-proximal  alternating coordinate method was proposed for solving such a structured convex-concave minimax problem and proved the global convergence as well as the linear rate of convergence.  Further, a primal-dual proximal splitting (PDPS) method proposed by Clason \emph{et al.} \cite{Clason2021} and an inexact primal-dual smoothing (IPDS) framework studied by Hien \emph{et al.} in \cite{Hien2017} are used to solve minimax problems with nonsmooth terms, and the convergence is established provided with some Lipschitz properties. Recently, Hamedani \cite{Hamedani2021} proposed a primal-dual algorithm for the nonsmooth convex-concave minimax  problem and achieved an ergodic convergence rate of function value with $O(1/k)$. However, these methods are not invalid  when $x$ and $y$ interact in constraints. For nonsmooth linearly nonseparable constraint minimax problems of the form (\ref{minimax}), even for the smooth problem (\ref{minimax}) with $\varphi\equiv 0$ and  $\psi \equiv 0$, to our best knowledge, numerical algorithms are quite rare. For smooth minimax problems, Tsaknakis \emph{et al.} \cite{Tsaknakis2021} proposed a multiplier gradient descent method and studied its convergence properties. However, they did not analyze iteration complexities of the method in \cite{Tsaknakis2021}.  Therefor, we hope to develop numerical algorithms for nonsmooth linearly nonseparable constraint minimax problems, establish their iteration complexities and apply them to solve linear equations with absolute values or projection operators.

\subsection{Contribution}
In this paper, we focus on optimality conditions and numerical algorithms for nonsmooth linearly nonseparable constraint convex-concave minimax problems \eqref{minimax}. One of obstacles is to deal with the nonseparable constraint in \eqref{minimax}. It is worth mentioning that under mild conditions (Assumption \ref{assum-1e0}, i.e., some Lipschitz properties), the above problem \eqref{minimax} can be reformulated as the following unconstrained nonsmooth minimax problem
$$
{\displaystyle \min_{x,\lambda}\max_y \varphi (x)+g(x)+x^TKy-h(y)-\psi (y)+\langle \lambda, Ax+By+c\rangle,}
$$
where $\lambda\in\mathfrak{R}^q$ is a Lagrange multiplier. For the above unconstrained nonsmooth problem, inspired by the success of proximal gradient methods in nonsmooth optimization, we use the proximal mapping technique to deal with nonsmooth terms and propose proximal gradient methods. The convergence of algorithms depends on finding an $\epsilon$-stationary point of the subproblem with respect to $y$.

The main contribution of this paper can be summarized as follows.
\begin{itemize}
  \item We define the stationary point of \eqref{minimax} and established equivalence conditions and properties of the stationary point;
  \item We develope an alternating coordinate ascent-decent algorithm framework and proved the convergence of iterations;
  \item We propose a proximal gradient multi-step ascent decent method and showed that this algorithm can find an $\epsilon$-stationary point in ${\cal O}\left(\epsilon^{-2}\log  \epsilon^{-1}  \right)$  iterations;
  \item We prove that GAVE and GLPE are equivalent to nonseparable linearly constraint convex-concave minimax problems and used the proposed numerical algorithm to solve these popular linear equations.
\end{itemize}

\subsection{Organization}
The paper is organized as follows. In Section 2, we provide some technique results about  minimizing the sum of  a nonsmooth  function and a ${\cal C}^{1,1}$ function, which play important roles in analyzing convergence properties of the proposed algorithms. In Section 3, we define the  stationary point of nonsmooth linearly constrained convex-concave
minimax problems and  establish optimality conditions.
In Section 4, we develop an alternating coordinate ascent-decent method, and any accumulation point of the sequences generated by this method is a stationary point of (\ref{minimax}). In Section 5, we construct  a proximal gradient multi-step ascent decent method for (\ref{minimax}) and prove that this algorithm generates a sequence converging to an $\epsilon$-stationary point within ${\cal O}(\epsilon^2 \log \epsilon^{-1})$ iterations. In Section 6, we apply the proximal gradient multi-step ascent decent
method to generalized absolute value equations and linear regression problems and show the effect of
our algorithm.  Some remarks are made in the last section.

\subsection{Notations}
We use the following notations throughout this paper. $\mathfrak{R}^n_+$ denotes $n$-dimensional positive octant cone. We use $\|\cdot\|$ to represent the Euclidean norm and its induced matrix norm and $\textbf{B}\left(x,d\right)$  to denote a closed ball of radius $d>0$ centered at $x $. For any $x\in\mathfrak{R}^n$, $\nabla h(x)$ and $\partial g(x)$ represent the gradient of the smooth function $h$ at $x$ and the sub-derivative of the nonsmooth function $g$ at $x$, respectively. The conjugate function of $f:\mathfrak{R}^n\rightarrow\mathfrak{R}^m$ is defined as $f^*:\mathfrak{R}^m\rightarrow\mathfrak{R}^n$.
\section{Preliminary}\label{Sec2}
\setcounter{equation}{0}
\quad \, In this section, we give several  results about properties of minimizing the sum of a nonsmooth  function and a ${\cal C}^{1,1}$ function. The results will be used in the following sections. Consider the extended real-valued function $F:\Re^d \rightarrow \overline \Re$ of the form
$$
F(z)=\sigma (z)+h(z),
$$
where $\sigma:\Re^d \rightarrow \overline \Re$ is a proper lower semicontinuous function and $h: \Re^d \rightarrow\Re$ is an $L_h$-smooth function, {\it i.e.}, $h$ is continuously differentiable  whose gradient mapping is Lipschitz continuous with Lipschitz constant $L_h>0$:
$$
\|\nabla h(z')-\nabla h(z)\|\leq L_h \|z'-z\|, \quad \forall z',z \in \Re^d.
$$

For $\lambda>0$, the proximal mapping for $\sigma$, denoted as ${\rm prox}_{\lambda \sigma}$, is defined by
$$
{\rm prox}_{\lambda \sigma}(z)={\rm argmin}\left\{\lambda\sigma (z')+\displaystyle \frac{1}{2}\|z'-z\|^2\right\}.
$$
Define, for $L>0$
$$
T_L^{h,\sigma}(z)={\rm prox}_{L^{-1} \sigma}(z-L^{-1} \nabla h(z))
$$
and
$$
G_{L}^{h,\sigma}(z)=L(z-T_L^{h,\sigma}(z)).
$$
We first state the descent lemma for $L_h$-smooth function, see Lemma 5.7 of \cite{Beck2017}.
\begin{lemma}\label{lem0}
Let $h$ be an $L_h$-smooth function defined over an open convex set $C \subseteq \Re^d$. Then for any $z,z'\in C$,
$$
h(z')\leq h(z)+\langle \nabla h(z),z'-z\rangle+\displaystyle \frac{L_h}{2}\|z'-z\|^2.
$$
\end{lemma}

The following results are taken from Section 10.3 of \cite{Beck2017}.
\begin{lemma}\label{lem1}
Let $\sigma:\Re^d \rightarrow \overline \Re$ is a proper lower semicontinuous convex function and $h: \Re^d \rightarrow\Re$ is a smooth function whose gradient mapping is Lipschitz continuous with Lipschitz constant $L_h>0$. Then, for $F=\sigma+h$ and $z \in {\rm dom} \, \sigma$,
\begin{itemize}
\item[{\rm (a)}] The point $z^*$ is a stationary point of $F$ if and only if $G_L^{h,\sigma}(z^*)=0$;
\item[{\rm (b)}] $\|G_L^{h,\sigma}(z')-G_L^{h,\sigma}(z)\|\leq (2L+L_h)\|z'-z\|$;
\item[{\rm (c)}]$F(z)-F\left(T^{h,\sigma}_{L}(z)\right)\geq \displaystyle \frac{2L-L_h}{2L^2}\|G_L^{h,\sigma}(z)\|^2$;
\item[{\rm (d)}] $F(z)-F\left(T^{h,\sigma}_{L_h}(z)\right)\geq \displaystyle \frac{1}{2L_h}\|G_{L_h}^{h,\sigma}(z)\|^2$;
\item[{\rm (e)}] If in addition $h$ is convex, then for $L \geq L_h$ and $z\in {\rm dom}\, \sigma$,
$$
\|G_{L}^{h,\sigma}(T_{L}^{h,\sigma}(z))\|\leq \|G_{L}^{h,\sigma}(z)\|.
$$
\end{itemize}
\end{lemma}

The following lemma is a modified version of Lemma 2 of Bolte \cite{Bolte2014}.
\begin{lemma}\label{lem2}
Let $\sigma:\Re^d \rightarrow \overline \Re$ is a proper lower semicontinuous function with
$$
\inf_{z \in \Re^d} \sigma (z)>-\infty,
$$
and $h: \Re^d \rightarrow\Re$ is a smooth function whose gradient mapping is Lipschitz continuous with Lipschitz constant $L_h>0$. Then, for $F=\sigma+h$, $z \in {\rm dom} \, \sigma$, $\xi \in \Re^d$ and
$$
z^+ \in {\rm prox}_{t^{-1}\sigma}\left(z-\displaystyle \frac{1}{t}\xi\right),
$$
we have
$$
F(z^+) \leq F(z)-\displaystyle \frac{1}{2}(t-L_h)\|z^+-z\|^2+\langle \nabla h(z)-\xi, z^+-z\rangle.
$$
\end{lemma}
{\bf Proof}. From the definition of $z^+$, one has that
$$
z^+ \in {\rm argmin}\left\{\langle \xi,z'-z\rangle+\displaystyle \frac{t}{2}\|z'-z\|^2+\sigma (z'): z'\in \Re^d \right\}.
$$
This implies, for any $z'\in \Re^d$, that
$$
\langle \xi,z'-z\rangle+\displaystyle \frac{t}{2}\|z'-z\|^2+\sigma (z')\geq \langle \xi,z^+-z\rangle+\displaystyle \frac{t}{2}\|z^+-z\|^2+\sigma (z^+).
$$
Taking $z'=z$ in the above inequality yields
$$
\sigma (z)\geq \langle \xi,z^+-z\rangle+\displaystyle \frac{t}{2}\|z^+-z\|^2+\sigma (z^+).
$$
Since $h$ is continuously differentiable and $\nabla h$ is $L_h$-Lipschitz continuous, we get
$$
h(z^+) \leq h(z)+\langle \nabla h(z),z^+-z \rangle+\displaystyle \frac{L_h}{2} \|z^+-z\|^2.
$$
Combining the above two inequalities, we obtain
$$
\sigma (z^+)+h(z^+)\leq \sigma (z)+h(z)-\displaystyle \frac{1}{2}(t-L_h)\|z^+-z\|^2+\langle \nabla h(z)-\xi, z^+-z\rangle.
$$
The proof is completed. \hfill $\Box$\\

The following lemma is a modified version of Theorem 10.16  of \cite{Beck2017}.
\begin{lemma}\label{lem2a}
Let $\sigma:\Re^d \rightarrow \overline \Re$ be a proper lower semicontinuous convex function
and $h: \Re^d \rightarrow\Re$ be a smooth function whose gradient mapping is Lipschitz continuous with Lipschitz constant $L_h>0$. Then, for $F=\sigma+h$, $z \in {\rm dom} \, \sigma$, $t \geq L_h$ and
$$
z^+ = {\rm prox}_{t^{-1}\sigma}\left(z-\displaystyle \frac{1}{t}\nabla h(z)\right),
$$
we have for any $z' \in {\rm dom}\,\sigma$,
$$
F(z')-F(z^+) \geq \displaystyle \frac{t}{2}\|z'-z^+\|^2-\displaystyle \frac{t}{2}\|z'-z\|^2+l_h(z',z),
$$
where
$$
l_h(z',z)=h(z')-h(z)-\langle \nabla h(z),z'-z\rangle.
$$
\end{lemma}
{\bf Proof}. Consider the function
$$
\phi (z')=h(z)+\langle \nabla h(z),z'-z\rangle +\sigma (z')+\displaystyle \frac{t}{2}\|z'-z\|^2.
$$
Since $\phi$ is a $t$-strongly convex function and $z^+={\rm argmin}_{z'} \phi(z')$, one has that
$$
\phi (z')-\phi (z^+) \geq \displaystyle \frac{t}{2}\|z'-z^+\|^2.
$$
Noting that
$$
\begin{array}{rcl}
\phi (z^+)& = &h(z)+\langle \nabla h(z),z^+-z\rangle +\sigma (z^+)+\displaystyle \frac{t}{2}\|z^+-z\|^2\\[6pt]
& \geq & h(z)+\langle \nabla h(z),z^+-z\rangle +\sigma (z^+)+\displaystyle \frac{L_h}{2}\|z^+-z\|^2\\[6pt]
& \geq & h(z^+)+\sigma (z^+)=F(z^+),
\end{array}
$$
we have for any $z'\in {\rm dom}\,\sigma$,
$$
\phi (z')-F(z^+) \geq \displaystyle \frac{t}{2}\|z'-z^+\|^2.
$$
Plugging the expression of $\phi (z')$ into the above inequality, we obtain
$$
h(z)+\langle \nabla h(z),z'-z\rangle+\sigma (z')+\displaystyle \frac{t}{2}\|z'-z\|^2-F(z^+)\geq
\displaystyle \frac{t}{2}\|z'-z^+\|^2,
$$
which is the same as the desired inequality
$$
F(z')-F(z^+) \geq  \displaystyle \frac{t}{2}\|z'-z^+\|^2-\displaystyle \frac{t}{2}\|z'-z\|^2+h(z')-h(z)-\langle \nabla h(z),z'-z\rangle.
$$
The proof is completed. \hfill $\Box$

\section{Properties of the nonsmooth minimax optimization problem}
\setcounter{equation}{0}
In this section, we analyze properties of the nonsmooth  minimax optimization problem (\ref{minimax}).
Define the feasible region of Problem (\ref{minimax}) by
$$
C=\{(x,y) \in \Re^n \times \Re^m: Ax+By+c=0\}.
$$
 First of all, we propose some assumptions about functions $g,h,\varphi$ and $\psi$.
 Denote $\varphi_g:=\varphi+g$ and $\psi_h:=h+\psi$.
\begin{assumption}\label{assum-1e0}
Let functions $g:\Re^n  \rightarrow \Re$,$h:\Re^m  \rightarrow \Re$, $\varphi:\Re^n \rightarrow \overline \Re$ and $\psi:\Re^m \rightarrow \overline \Re$ satisfy the following conditions.
\begin{itemize}
\item[{\bf A1}] Functions $g$ and $h$ are continuously differentiable convex with Lipschitz continuous gradients, {\it i.e.}, there exist constants $L_g>0$ and $L_h>0$ such that
$$
\|\nabla g(x')-\nabla g(x)\|\leq L_g \|x'-x\|,\, \|\nabla h(y')-\nabla h(y)\|\leq L_h \|y'-y\|,\,
 \forall (x',y'), (x,y) \in \Re^n \times \Re^m.
$$
\item[{\bf A2}] Functions $\varphi$ and $\psi$ are proper lower semicontinuous convex functions.
\item[{\bf A3}] Function $\psi_h$ is a $\mu$-strongly convex function for some positive constant $\mu>0$.
\end{itemize}
\end{assumption}

Define the Lagrange function of \eqref{minimax}
\begin{equation*}
  \mathcal {L}(x,y,\lambda)=\varphi(x)+g(x)+x^TKy-h(y)-\psi (y)+\langle \lambda, Ax+By+c\rangle,
\end{equation*}
where $\lambda\in \mathfrak{R}^q$ is a Lagrange multiplier. For a fixed point $(x,\lambda)$, we define functions with respect to  $y$ as follows
\begin{equation}\label{eq:Qtheta}
\begin{array}{l}
Q(x,y,\lambda)=x^TKy-h(y)-\psi (y)+\langle \lambda, Ax+By+c\rangle,\\[6pt]
 \theta_0(x,\lambda)=\max_y Q(x,y,\lambda),\, y_*(x,\lambda)={\rm argmax}_y Q(x,y,\lambda),\\[6pt]
\theta (x,\lambda)=g(x)+\theta_0(x,\lambda).
\end{array}
\end{equation}
 If Assumption \ref{assum-1e0} is satisfied, then from the definition of $y_*(x,\lambda)$, we have
\begin{equation}\label{eq:h01}
0\in -[K^Tx+B^T\lambda]+\partial \psi_h (y_*(x,\lambda)).
\end{equation}

\begin{remark}\label{remarkCD}
It follows from classical convex analysis that  Assumption \ref{assum-1e0} implies that the conjugate function $\psi_h^*$ is continuously differentiable with $\mu^{-1}$-Lipschitz gradient.
From the Moreau-Fenchel theorem \cite[Proposition 11.3]{RW98} for a proper closed function,
(\ref{eq:h01})  implies that
$$
y_*(x,\lambda)\in \partial \psi_h^*(K^Tx+B^T\lambda),
$$
which with  the differentiability of $\psi_h^*$ gives
\begin{equation}\label{eq:ystar}
y_*(x,\lambda)=\nabla \psi_h^*(K^Tx+B^T\lambda).
\end{equation}
\end{remark}

 Now we define
\begin{equation}\label{eq:gm}
\gamma=\max\left\{\sqrt{2\left(L_g\mu+\|K\|^2\right)^2+2\left(\|A\|\mu+\|K\|\|B\|\right)^2},\sqrt{2{(\mu\|A\|+\|K\|\|B\|)^2+2\|B\|^4}}\right\}.
\end{equation}
In the following proposition, we will prove that, under  Assumption \ref{assum-1e0}, $y_*$ and $\nabla \theta$ are both Lipschitz continuous with Lipschitz constants $\displaystyle \frac{1}{\mu}$ and  $\displaystyle \frac{\gamma}{\mu}$, respectively.
\begin{proposition}\label{prop-1e0}
Let Assumption  \ref{assum-1e0} be satisfied. One has, for any
$(x^0,\lambda^0)\in {\rm dom}\, \theta$ and $(x^1,\lambda^1)\in {\rm dom}\, \theta$, that
\begin{equation}\label{eq:Ae0}
\|y_*(x^0,\lambda^0)-y^*(x^1,\lambda^1)\|\leq \displaystyle \frac{1}{\mu} [\|K\|\|x^0-x^1\|+\|B\|\|\lambda^0-\lambda^1\|].
\end{equation}
The function $\theta_0 (x,\lambda)$ is continuously differentiable at any $(x,\lambda) \in {\rm dom}\,\theta_0$ with
\begin{equation}\label{eq:gTheta0}
\nabla \theta (x,\lambda)=\left
(
\begin{array}{l}
\nabla g(x)+A^T\lambda+K\nabla \psi_h^*(K^Tx+B^T\lambda)\\
Ax+c+B\nabla \psi_h^*(K^Tx+B^T\lambda)
\end{array}
\right).
\end{equation}
and
\begin{equation}\label{eq:Be0}
\|\nabla \theta(x^0,\lambda^0)-\nabla \theta(x^1,\lambda^1)\|\leq \displaystyle \frac{\gamma}{\mu}\|(x^0,z^0)-(x^1,z^1)\|
\end{equation}
for any
$(x^0,z^0)\in {\rm dom}\, \theta$ and $(x^1,z^1)\in {\rm dom}\, \theta$.
\end{proposition}
{\bf Proof}. For simplicity, we denote $y^0=y_*(x^0,\lambda^0)$ and $y^1=y_*(x^1,\lambda^1)$. From the definition of $y_*(x,\lambda)$ and the $\mu$-strong concavity of the function $y \rightarrow Q(x,y,\lambda)$, we have for any $y\in \Re^m$,
\begin{equation}\label{eq-1e0}
-x^{0T}Ky+\psi_h(y)-\langle \lambda^0,Ax^0+By+c\rangle \geq -x^{0T}Ky^0+\psi_h(y^0)-\langle \lambda^0,Ax^0+By^0+c\rangle +\displaystyle
\frac{\mu}{2}\|y-y^0\|^2
\end{equation}
and
\begin{equation}\label{eq-2e0}
-x^{1T}Ky+\psi_h(y)-\langle \lambda^1,Ax^1+By+c\rangle \geq -x^{0T}Ky^1+\psi_h(y^1)-\langle \lambda^1,Ax^1+By^1+c\rangle +\displaystyle
\frac{\mu}{2}\|y-y^1\|^2.
\end{equation}
Replacing $y=y^1$ in (\ref{eq-1e0}) yields
\begin{equation}\label{eq-3e0}
-x^{0T}Ky^1+\psi_h(y^1)-\langle \lambda^0,Ax^0+By^1+c\rangle \geq -x^{0T}Ky^0+\psi_h(y^0)-\langle \lambda^0,Ax^0+By^0+c\rangle +\displaystyle
\frac{\mu}{2}\|y^1-y^0\|^2.
\end{equation}
Replacing $y=y^0$ in (\ref{eq-2e0}) yields
\begin{equation}\label{eq-4e0}
-x^{1T}Ky^0+\psi_h(y^0)-\langle \lambda^1,Ax^1+By^0+c\rangle \geq -x^{0T}Ky^1+\psi_h(y^1)-\langle \lambda^1,Ax^1+By^1+c\rangle +\displaystyle
\frac{\mu}{2}\|y^0-y^1\|^2.
\end{equation}
Adding both sides of (\ref{eq-3e0}) and  (\ref{eq-4e0}), we obtain
$$
\begin{array}{l}
\langle K^Tx^0,y^0\rangle-\langle K^Tx^0,y^1\rangle+\langle K^Tx^1,y^1\rangle-\langle K^Tx^1,y^0\rangle
+\langle \lambda^1-\lambda^0,B(y^1-y^0)\rangle\geq \mu \|y^0-y^1\|^2,
\end{array}
$$
or equivalently
\begin{equation}\label{eq-6e0}
\langle  K(y^1-y^0),x^1-x^0\rangle +\langle \lambda^1-\lambda^0,B(y^1-y^0)\rangle\geq \mu \|y^0-y^1\|^2.
\end{equation}
Therefore, we have from (\ref{eq-6e0})  that
$$
\mu \|y^0-y^1\|^2 \leq  \|K\|\|x^0-x^1\|\|y^0-y^1\|+\|B\| \|\lambda^0-\lambda^1\|\|y^0-y^1\|,
$$
which implies the inequality (\ref{eq:Ae0}).
It follows from Remark \ref{remarkCD} that $\psi_h^*$ is differentiable with $\mu^{-1}$-Lipschitz continuous gradient. Thus
$\theta$ is differentiable at $(x,\lambda) \in {\rm dom}\,\theta$ with
$$
\nabla \theta (x,\lambda)=\left
(
\begin{array}{l}
\nabla g(x)+A^T\lambda+K\nabla \psi_h^*(K^Tx+B^T\lambda)\\
Ax+c+B\nabla \psi_h^*(K^Tx+B^T\lambda)
\end{array}
\right).
$$
From this formula and the $\mu^{-1}$-Lipschitz continuity of $\nabla \psi_h^*$ , we have the following estimates
$$
\begin{array}{ll}
&\!\!\!\!\|\nabla \theta (x^0,\lambda^0)-\nabla \theta (x^1,\lambda^1)\|^2 \\[4pt]
&=\|\nabla g(x^0)-\nabla g(x^1){+A^T(\lambda^0-\lambda^1)}+K[
\nabla \psi_h^*(K^Tx^0+B^T\lambda^0)-\nabla \psi_h^*(K^Tx^1+B^T\lambda^1]\|^2\\[4pt]
& \quad +\|A(x^0-x^1)+B[
\nabla \psi_h^*(K^Tx^0+B^T\lambda^0)-\nabla \psi_h^*(K^Tx^1+B^T\lambda^1]\|^2\\[4pt]
&\leq [L_g\|x^1-x^0\|+{\|A\|\|\lambda^0-\lambda^1\|}+\mu^{-1}\|K\|(\|K\|\|x^1-x^0\|+\|B\|\lambda^0-\lambda^1\|)]^2\\[4pt]
& \quad + [\|A\|\|x^1-x^0\|+\mu^{-1}\|B\|(\|K\|\|x^1-x^0\|+\|B\|\lambda^0-\lambda^1\|)]^2\\[4pt]
&= [(L_g+\mu^{-1}\|K\|^2)\|x^1-x^0\|+{(\|A\|+\mu^{-1}\|K\|\|B\|)}\|\lambda^0-\lambda^1\|]^2\\[4pt]
&\quad +
[(\|A\|+\mu^{-1}\|K\|\|B\|)\|x^1-x^0\|+\mu^{-1}\|B\|^2\|\lambda^0-\lambda^1\|]^2\\[4pt]
& \leq \displaystyle \frac{2}{\mu^2}\left[(L_g\mu+\|K\|^2)^2+(\|A\|\mu+\|K\|\|B\|)^2\right]\|x^1-x^0\|^2\\[6pt]
& \quad +\displaystyle \frac{2}{\mu^2}\left[{(\mu\|A\|+\|K\|\|B\|)^2+\|B\|^4}\right]\|\lambda^1-\lambda^0\|^2\\[10pt]
& \leq  \displaystyle \frac{\gamma^2}{\mu^2}[\|x^1-x^0\|^2+\|\lambda^1-\lambda^0\|^2],
\end{array}
$$
where $\gamma$ is defined by (\ref{eq:gm}). This proves the estimate (\ref{eq:Be0}).
The proof is completed. \hfill $\Box$\\

Define
\begin{equation}\label{eq:thtaL}
L_{\theta}= \displaystyle \frac{\gamma}{\mu}.
\end{equation}
Then, $\theta$ is continuously differentiable and $\nabla \theta$ is $L_\theta$-Lipschitz continuous when  Assumption  \ref{assum-1e0} is satisfied. Next, we establish the optimality conditions of Problem (\ref{minimax}). We define the stationary point of Problem (\ref{minimax}), which is similar to KKT point for constrained optimization problems.
\begin{definition}\label{spointe0}
We say that $(\overline x,\overline y)$ is a stationary point of Problem (\ref{minimax})  if there exists a vector  $\overline \lambda \in \Re^q$ such that
\begin{equation}\label{eq:NCs}
0\in  K\overline y+A^T\overline \lambda +\nabla g(\overline x)+\partial \varphi(\overline x),\,\,\, 0\in K^T\overline x+B^T\overline \lambda-\nabla h(\overline y)-\partial \psi(\overline y),\,\,\, A\overline x+B\overline y+c=0.
\end{equation}
\end{definition}

\begin{proposition}\label{prop-2e0}
Let Assumption  \ref{assum-1e0} be satisfied. Then
$$0\in \nabla \theta (\overline x,\overline \lambda)+\partial \varphi (\overline x)\times \{0\}$$ if and only if
$\left(\overline x,\overline y=\nabla {\psi}_h^*\left(A^T\overline x+B^T \overline \lambda\right)\right)$ is a stationary point of Problem (\ref{minimax}).
\end{proposition}
{\bf Proof}. From Proposition \ref{prop-1e0} and the relation (\ref{eq:ystar}), we have
\begin{equation}\label{eq:h02}
\begin{array}{ll}
\nabla \theta (x,\lambda)& =\left
(
\begin{array}{l}
\nabla g(x)+A^T\lambda+K\nabla \psi_h^*(K^Tx+B^T\lambda)\\
Ax+c+B\nabla \psi_h^*(K^Tx+B^T\lambda)
\end{array}
\right)\\[10pt]
&=\left
(
\begin{array}{l}
\nabla g(x)+A^T\lambda+Ky_*(x,\lambda)\\
Ax+c+By_*(x,\lambda)
\end{array}
\right).
\end{array}
\end{equation}
 If $ 0\in \nabla \theta (\overline x,\overline \lambda)+\partial \varphi (\overline x)\times \{0\}$, we have from (\ref{eq:h02}) and (\ref{eq:h01}) that
$$
\begin{array}{l}
0\in -[K^T\overline x+B^T\overline \lambda]+\partial \psi_h (y_*(\overline x,\overline \lambda)),\\[6pt]
0\in \nabla g(\overline x)+A^T\overline \lambda+Ky_*(\overline x,\overline \lambda)+\partial \varphi (\overline x),\\[6pt]
0=A\overline x+c+By_*(\overline x,\overline \lambda)
\end{array}
$$
Therefore, for $\overline y=\nabla \psi_h^*\left(A^T\overline x+B^T \overline \lambda\right)$, we have $\overline y=y_*(\overline x,\overline \lambda)$ and $\left(\overline x,\overline y=\nabla \psi_h^*\left(A^T\overline x+B^T \overline \lambda\right)\right)$ is a stationary point of Problem (\ref{minimax}).

On the other hand, let $(\overline x,\overline y,\overline \lambda)$ satisfy (\ref{eq:NCs}), or equivalently
\begin{equation}\label{eq:h0x}
 0\in  K\overline y+A^T\overline \lambda +\partial \varphi_g(\overline x),\,\,\, 0\in K^T\overline x+B^T\overline \lambda-\partial \psi_h(\overline y),\,\,\, A\overline x+B\overline y+c=0.
\end{equation}
Then, from $0\in K^T\overline x+B^T\overline \lambda-\partial \psi_h(\overline y)$, we have from the Moreau-Fenchel equality for a proper closed function that
$$
\overline y=\nabla \psi_h^*(K^T\overline x+B^T\overline \lambda)
$$
and $\nabla \psi_h^*(K^T\overline x+B^T\overline \lambda)=y_*(\overline x,\overline \lambda)$.
In this case,
$$
\nabla \theta (\overline x,\overline \lambda)=\left
(
\begin{array}{l}
\nabla g(\overline x)+A^T\overline \lambda+K\nabla \psi_h^*(K^T\overline x+B^T\overline\lambda)\\
A\overline x+c+B\nabla \psi_h^*(K^T\overline x+B^T\overline \lambda)
\end{array}
\right),
$$
which implies $0\in \nabla \theta (\overline x,\overline \lambda)+\partial \varphi (\overline x)\times \{0\}$ from (\ref{eq:h0x}).
The proof is completed. \hfill $\Box$

In view of  (\ref{eq:Be0}), we have that
$$
\theta(x,\lambda)=g(x)+c^T\lambda+\langle \lambda, Ax\rangle+\psi_h^*(K^Tx+B^T\lambda)
$$
is continuously differentiable with $L_{\theta}$-Lipschitz continuous gradient mapping.
\begin{proposition}\label{prop:thetaJ}
Let Assumption  \ref{assum-1e0} be satisfied. Then for any $(x,\lambda) \in {\rm dom}\, \theta$,
\begin{equation}\label{eq:GJ}
\partial \nabla \theta (x,\lambda)=\left\{ \left[
\begin{array}{ll}
\nabla^2 g(x) & A^T\\[4pt]
A & 0
\end{array}
\right]+\left[
\begin{array}{l}
K\\[4pt]
B
\end{array}
\right]V[K^T \quad B^T]: V \in \partial \nabla \psi_k^*(K^Tx+B^T\lambda)\right\},
\end{equation}
where $\partial$ denotes the generalized Jacobian in the sense of Clarke.
\end{proposition}
{\bf Proof}. Let ${\cal D}_{\nabla\psi_h^*}$ and ${\cal D}_{\nabla \theta}$ denote the sets of differentiable points of locally Lipschitz continuous functions $\nabla \psi_h^*$ and $\nabla \theta$, respectively. Then we have
 $$
 {\cal D}_{\nabla \theta}=\{(x,\lambda):K^Tx+B^T\lambda \in {\cal D}_{\nabla\psi_h^*}\}.
 $$
  For any $(x,\lambda) \in {\cal D}_{\nabla\theta}$, the Jacobian of $\theta$ at $(x,\lambda)$ is expressed as
$$
\begin{array}{ll}
\nabla^2 \theta (x,\lambda)& =\left[
\begin{array}{ll}
\nabla^2 g(x)+K \nabla^2 \psi_h^*(K^Tx+B^T\lambda)K^T & A^T+K \nabla^2 \psi_h^*(K^Tx+B^T\lambda)B^T\\[4pt]
A+B \nabla^2 \psi_h^*(K^Tx+B^T\lambda)K^T & B\nabla^2 \psi_h^*(K^Tx+B^T\lambda)B^T
\end{array}
\right]\\[16pt]
&=\left[
\begin{array}{ll}
\nabla^2 g(x) & A^T \\[4pt]
A & 0
\end{array}
\right]+\left[
\begin{array}{l}
K\\[4pt]
B
\end{array}
\right]
\nabla^2 \psi_h^*(K^Tx+B^T\lambda)
[K^T \quad B^T].
\end{array}
$$
From this expression, we can easily obtain\footnote{The following $\limsup$ denotes the outer limit of a set-valued mapping in the sense of Rockafellar and Wets (1998) \cite{RW98}.}
$$
\begin{array}{ll}
\partial_{B} \nabla \theta (x,\lambda) &= \displaystyle\limsup_{(x',\lambda')  \stackrel{{\cal D}_{\nabla\theta}}\rightarrow  (x,\lambda)}\nabla \theta (x',\lambda')\\[18pt]
&=\left\{\left[
\begin{array}{ll}
\nabla^2 g(x) & A^T \\[4pt]
A & 0
\end{array}
\right]+\left[
\begin{array}{l}
K\\[4pt]
B
\end{array}
\right]
V
[K^T \quad B^T]:V\in \partial_B\nabla \psi_h^*(K^Tx+B^T\lambda)  \right\}.
\end{array}
$$
Noting $\partial \nabla \theta (x,\lambda)={\rm co}\left(\partial \nabla \theta (x,\lambda) \right)$, we obtain formula
(\ref{eq:GJ}).  The proof is completed. \hfill $\Box$

Let
\begin{equation}\label{eq:fdefi}
f(x,y,\lambda)=g(x)+x^TKy-h(y)+\langle \lambda, Ax+By+c\rangle.
\end{equation}
Then, under Assumption \ref{assum-1e0},
Problem (\ref{minimax}) is reformulated as the following unconstrained nonsmooth minimax optimization problem
\begin{equation}\label{prob-RP}
{\displaystyle \min_{x,\lambda}\max_y \mathcal {L}(x,y,\lambda)=\varphi (x)+f(x,y,\lambda)-\psi (y).}
\end{equation}
Obviously, we have that $(\overline x, \overline y,\overline \lambda)$ is a stationary point of Problem (\ref{prob-RP}) if and only if $(\overline x, \overline y,\overline \lambda)$ satisfies the following system
$$
\left
\{
\begin{array}{l}
0 \in \nabla g(\overline x)+K \overline y+A^T \overline \lambda+\partial \varphi (\overline x),\\[4pt]
0=A\overline x+B\overline y+c,\\[4pt]
0\in \nabla h(\overline y)-K^T \overline x-B^T\overline \lambda +\partial \psi (\overline y),
\end{array}
\right.
$$
which coincides with the set of conditions (\ref{eq:NCs}). As $\partial \varphi$ and $\partial \psi$ are hard to be portrayed, proximal mapping is used to characterize the properties of the stationary point of Problem (\ref{minimax}). Hence, let $L>0$ and define
$$
\begin{array}{l}
G^{f,\varphi}_L(x,y,\lambda)=L(x-{\rm prox}_{L^{-1}\varphi}(x-L^{-1}\nabla_x f(x,y,\lambda))),\\[6pt]
{G^{f,0}_L(x,y,\lambda)=L(\lambda-{\rm prox}_{L^{-1}0}(\lambda-L^{-1}\nabla_\lambda f(x,y,\lambda))),}\\[6pt]
G^{f,\psi}_L(x,y,\lambda)=L(y-{\rm prox}_{L^{-1}\psi}(y+L^{-1}\nabla_yf(x,y,\lambda))).
\end{array}
$$
Then we have
$$
\begin{array}{l}
G^{f,\varphi}_L(x,y,\lambda)=L\left(x-{\rm prox}_{L^{-1}\varphi}\left(x-L^{-1}\left(\nabla g(x)+Ky+A^T\lambda\right)\right)\right),\\[6pt]
G^{f,0}_L(x,y,\lambda)=Ax+By+c,\\[6pt]
G^{f,\psi}_L(x,y,\lambda)=L\left(y-{\rm prox}_{L^{-1}\psi}\left(y+L^{-1}\left(-\nabla h(y)+K^Tx+B^T\lambda\right)\right)\right).
\end{array}
$$
It is not hard to verify the following reslut.
\begin{lemma}\label{lem:optimality}
 For $L_1>0$,$L_2>0$ and $L_3>0$, $(\overline x, \overline y)$ is a stationary point of Problem (\ref{minimax}) if and only if there exists $\overline\lambda \in \Re^q$ such that
$$
G^{f,\varphi}_{L_1}(\overline x, \overline y,\overline \lambda)=0,\ G^{f,\psi}_{L_2}(\overline x, \overline y,\overline \lambda)=0 \mbox{  and  }G^{f,0}_{L_3}(\overline x, \overline y,\overline \lambda)=0.
$$
\end{lemma}
\begin{definition}\label{aspointee0}
For $\epsilon>0$, we say that $(x,y)$ is an $\epsilon$-stationary point of Problem (\ref{minimax}) if and only if there exists a vector $\lambda\in \Re^q$ such that
$$
\|G^{f,\varphi}_{L_1}(x, y,\lambda)\|\leq \epsilon,\,\, \|G^{f,\psi}_{L_2}(x,y,\lambda)\|\leq \epsilon  \mbox{  and  } Ax+By+c=0
$$
for some positive constants $L_1$ and $L_2$.
\end{definition}

From the definitions of $f$ and $Q$, we have the following relation
$$
G^{f,\psi}_L(x,y,\lambda)=G^{Q,\psi}_L(x,y,\lambda)=L(y-{\rm prox}_{L^{-1}\psi}(y+L^{-1}\nabla_yQ(x,y,\lambda))).
$$
Since $y \rightarrow Q(x,y,\lambda)$ is $\mu$-strongly concave, we have, for any $L>0$, that
$$
G^{f,\psi}_L(x,y,\lambda)=G^{Q,\psi}_L(x,y,\lambda)=0
$$
if and only if $y=y_*(x,\lambda)$.
\section{An alternating coordinate method }
\setcounter{equation}{0}
The following algorithm is trying to get an approximate stationary point of Problem (\ref{minimax}), whose main idea is to construct an  alternating coordinate ascent-decent method for solving Problem (\ref{prob-RP}).
\begin{algorithm}\label{alg-30exact}
Input $(x^0,y^0,\lambda^0) \in \Re^n \times \Re^q \times \Re^m$, $\alpha_x >0$, $\varepsilon_t >0$ for $t \in \textbf{N}$\\
\mbox{}\hspace{0.2cm} {\bf for} $t=0,1,2,\ldots,$ {\bf do}\\[4pt]
\mbox{}\hspace{0.8cm}Find $y^{t+1}$ such that
\begin{equation}\label{yG}
\|G^{Q,\psi}_1(x^t,y^{t+1},\lambda^t)\|\leq \varepsilon_t;
\end{equation}
\mbox{}\hspace{0.8cm}Calculate
$$
\begin{array}{rcl}
x^{t+1}& =& {\rm prox}_{\alpha_x \varphi}\left[x^t-\alpha_x \left(\nabla g(x^t)+Ky^{t+1}+A^T\lambda^t\right)\right]\\[4pt]
\lambda^{t+1}& =&{\lambda^t-\alpha_x [Ax^t+By^{t+1}+c]}
\end{array}
$$
\mbox{}\hspace{0.2cm} {\bf end for}\\
\mbox{}\hspace{0.2cm} Return $(x^{t+1},y^{t+1},\lambda^{t+1})$ for $t=0,1,2,\ldots$
\end{algorithm}

Algorithm \ref{alg-30exact} is more like an algorithmic framework, where the technique of solving $y$ is not specified in \eqref{yG}. The convergence of Algorithm \ref{alg-30exact} only depends on whether the subproblem can be solved accurately enough. For simplicity, we denote $y_*(t)=y_*(x^t,\lambda^t)$. Then
$
G^{f,\psi}_L(x^t,y,\lambda^t)=0
$ for some constant $L>0$ if and only if $y=y_*(t)$.

Let $$\theta_{\varphi} (x,\lambda):=\varphi(x)+\theta (x,\lambda).$$
Then
$$
\begin{array}{rcl}
\theta_{\varphi} (x,\lambda) & = & \varphi(x)+g(x)+\displaystyle \sup_y \left[x^TKy-h(y)-\psi(y)+\langle \lambda,Ax+By+c\rangle\right]\\[6pt]
&=& \varphi_g(x)+\lambda^T(Ax+c)+\psi_h^*(K^Tx+B^T\lambda).
\end{array}
$$
To achieve the convergence properties of Algorithm \ref{alg-30exact}, we need the following assumptions.
\begin{assumption}\label{assum-2e0}
Suppose that $\varphi_g$, $\psi_h$ and $A,\ B,\ c$ satisfy
$$
\inf_{x,\lambda}\theta_{\varphi}(x,\lambda) >-\infty.
$$
\end{assumption}

For studying Problem (\ref{prob-RP}), Assumption \ref{assum-2e0} is reasonable when we try to find a local minimax point of unconstrained minimax optimization problem \eqref{prob-RP} in the sense of Jin, Netrapalli and Jordan \cite{Jin2019}.
\begin{assumption}\label{assum-eb}
Suppose that there exists a constant $\rho_0>0$ such that the following error bound condition holds
$$
\|y-y_*(t)\|\leq \rho_0 \|G^{f,\psi}_1(x^t,y,\lambda^t)\|, \quad \,\, \forall y \in \textbf{B}\left(y_*(t),\|y^t-y_*(t)\|\right).
$$
\end{assumption}

Assumption \ref{assum-eb} is a conventional error bound condition for the problem $\max_y Q(x^t,y,\lambda^t)$. In fact, Assumption \ref{assum-eb} holds trivially if A3 of Assumption \ref{assum-1e0} is satisfied and $\psi (y)\equiv 0$.\\
Let
\begin{equation}\label{eq:xit}
\xi^t=\left
[
\begin{array}{l}
\nabla_x f(x^t,y^{t+1},\lambda^t)\\[4pt]
\nabla_{\lambda}f(x^t,y^{t+1},\lambda^t)
\end{array}
\right]=\left[
\begin{array}{l}
\nabla g(x^t)+Ky^{t+1}+A^T\lambda^t\\[4pt]
Ax^t+By^{t+1}+c
\end{array}
\right].
\end{equation}
\begin{proposition}\label{prop-4exact}
Let Assumptions \ref{assum-1e0} and \ref{assum-eb} be satisfied.
Let $\alpha_x \in \left(0,1/L_{\theta}\right)$  and the sequence $\left\{(x^t,\lambda^t,y^t): t=0,1,\ldots,\right\}$ be generated by Algorithm \ref{alg-30exact}. Suppose $\|y^{t+1}-y_*(t)\|\leq \|y^t-y_*(t)\|$   and
\begin{equation}\label{eq:epsilonC}
\displaystyle \sum_{t=0}^{\infty} \varepsilon_t^2< \infty.
\end{equation}
Then for $t=0,1,\ldots$
\begin{itemize}
\item[{\rm (i)}]The following series is convergent
\begin{equation}\label{eq:ccd}
\displaystyle \sum_{t=0}^{\infty}\left\|(x^{t+1},\lambda^{t+1})-(x^{t},\lambda^{t})\right\|^2< \infty;
\end{equation}
\item[{\rm (ii)}]For $\sigma (x,\lambda)=\varphi (x)$, one has
\begin{equation}\label{eq:pmc}
G^{\theta,\sigma}_{\alpha_x^{-1}}(x^t,\lambda^t) \rightarrow 0, \, G^{Q,\psi}_{1}(x^t,y^t \lambda^t)\rightarrow 0;
\end{equation}
 \item[{\rm (iii)}]Let $(\overline x,\overline y,\overline \lambda)$ be any accumulation point of $\{(x^t,y^t,\lambda^t)\}$. Then $(\overline x,\overline y)$ is  a stationary point of Problem (\ref{minimax}) with Lagrange multiplier $\overline \lambda$.
\end{itemize}
\end{proposition}
{\bf Proof}. For $\xi^t$ defined by (\ref{eq:xit}), $(x^{t+1},\lambda^{t+1})$ can be expressed as
\begin{equation}\label{eq:xtht1}
(x^{t+1},\lambda^{t+1})={\rm prox}_{\alpha_x \sigma}((x^t,\lambda^t)-\alpha_x \xi^t),
\end{equation}
where $\sigma (x,\lambda)=\varphi (x)$.  Applying Lemma \ref{lem2} with $h=\theta$ and this $\sigma$, we obtain
$$
\begin{array}{rcl}
\varphi (x^{t+1})+\theta (x^{t+1},\lambda^{t+1}) & \leq & \varphi (x^t)+\theta (x^t,\lambda^t)-\displaystyle \frac{1}{2}(\alpha_x^{-1}-L_{\theta})\|(x^{t+1},\lambda^{t+1})- (x^t,\lambda^t)\|^2\\[8pt]
& &  +\langle \nabla \theta (x^t,\lambda^t)-\xi^t, (x^{t+1},\lambda^{t+1})- (x^t,\lambda^t)\rangle.
\end{array}
$$
Using the simple inequality $\langle a, b\rangle\leq \|a\|^2/2+\|b\|^2/2$ for vectors $a$ and $b$, we obtain
\begin{equation}\label{eq:hlp1}
\begin{array}{r}
\varphi(x^{t+1})+\theta (x^{t+1},\lambda^{t+1}) \leq \varphi(x^t)+\theta (x^t,\lambda^t)-\displaystyle \frac{1}{2}(\alpha_x^{-1}-L_{\theta})\|(x^{t+1},\lambda^{t+1})- (x^t,\lambda^t)\|^2\\[6pt]
\quad +\displaystyle \frac{1}{4}(\alpha_x^{-1}-L_{\theta})\|(x^{t+1},\lambda^{t+1})- (x^t,\lambda^t)\|^2
+\displaystyle \frac{1}{(\alpha_x^{-1}-L_{\theta})}
\|\nabla \theta (x^t,\lambda^t)-\xi^t\|^2.
\end{array}
\end{equation}
Noting that $\nabla_x Q(x^t,\cdot,\lambda^t)$ and $\nabla_{\lambda} Q(x^t,\cdot,\lambda^t)
$ are $\|K\|$-Lipschitz continuous and $\|B\|$-Lipschitz continuous, respectively, we obtain
$$
\|\nabla_x Q(x^t,y^{t+1},\lambda^t)-\nabla_x \theta_0 (x^t,\lambda^t)\|^2=\|\nabla_x Q(x^t,y^{t+1},\lambda^t)-\nabla_x Q(x^t,y_*(t),\lambda^t)\|^2\leq \|K\|^2
\|y^{t+1}-y_*(t)\|^2
$$
and
$$
\|\nabla_{\lambda} Q(x^t,y^{t+1},\lambda^t)-\nabla_{\lambda} \theta (x^t,\lambda^t)\|^2=\|\nabla_z Q(x^t,y^{t+1},\lambda^t)-\nabla_z Q(x^t,y_*(t),\lambda^t)\|^2\leq \|B\|^2
\|y^{t+1}-y_*(t)\|^2.
$$
Then, noting $\xi^t=\nabla_{x,\lambda}f(x^t,y^{t+1},\lambda^t)=\nabla g(x^t)+\nabla_x Q(x^t,{y^{t+1}},\lambda^t)$ and $\nabla \theta (x^t,\lambda^t)=\nabla_{x,\lambda}f(x^t,y_*(t),\lambda^t)$$=\nabla g(x^t)+\nabla_{x,\lambda}Q(x^t,y_*(t),\lambda^t)$, we get that
$$
\|\xi^t-\nabla \theta (x^t,\lambda^t)\|^2\leq [\|K\|^2+\|B\|^2]\|y^{t+1}-y_*(t)\|^2.
$$
Since $\|y^{t+1}-y_*(t)\|\leq \|y^t-y_*(t)\|$, from Assumption \ref{assum-eb}, this inequality implies
\begin{equation}\label{eq:xita2}
\|\xi^t-\nabla \theta (x^t,\lambda^t)\|^2\leq \rho_0^2  [\|K\|^2+\|B\|^2]\|G^{f,\psi}_1(x^t,y^{t+1},\lambda^t)\|^2
\leq \rho_0^2  [\|K\|^2+\|B\|^2]\varepsilon_t^2,
\end{equation}
where the last inequality comes from the definition of $y^{t+1}$.  Then, from (\ref{eq:hlp1}), we obtain the following inequality
\begin{equation}\label{eq:keyIneqExact}
\begin{array}{rcl}
\displaystyle \frac{1-\alpha_x L_{\theta}}{4\alpha_x}\|(x^{t+1},\lambda^{t+1})- (x^t,\lambda^t)\|^2& \leq &
\theta (x^t,\lambda^t)+\varphi (x^t)-\varphi(x^{t+1})+\theta (x^{t+1},\lambda^{t+1})\\[10pt]
& & +\displaystyle \frac{\rho_0^2 }{\alpha_x^{-1}-L_{\theta}} \left[\|K\|^2+\|B\|^2\right]\varepsilon_t^2
\end{array}
\end{equation}
For $T>1$,  summing (\ref{eq:keyIneqExact})
 over $t=0,1,\ldots, T-1$,
  we obtain from Assumption \ref{assum-2e0} that
 $$
\begin{array}{rcl}
\displaystyle \sum_{t=0}^{T-1}\|(x^{t+1},\lambda^{t+1})-(x^{t},\lambda^{t})\|^2 &
 \leq & \displaystyle \frac{4\alpha_x}{1-\alpha_xL_{\theta}}[\theta_{\varphi} (x^0,\lambda^0)-\theta_{\varphi} (x^T,\lambda^T)] +\displaystyle \frac{4\alpha_x^2\rho_0^2(\|K\|^2+\|B\|^2)}{(1-\alpha_xL_{\theta})^2}\sum_{t=0}^{T-1}\varepsilon_t^2\\[8pt]
  &\leq & \displaystyle \frac{4\alpha_x}{1-\alpha_xL_{\theta}}[\theta_{\varphi} (x^0,\lambda^0)-\inf \theta_\varphi]+\displaystyle \frac{4\alpha_x^2\rho_0^2(\|K\|^2+\|B\|^2)}{(1-\alpha_xL_{\theta})^2}\sum_{t=0}^{T-1}\varepsilon_t^2.
\end{array}
$$
 Therefore we obtain (\ref{eq:ccd}) of  (i)  from condition (\ref{eq:epsilonC}).

Now we turn to the proof of (ii). From (\ref{eq:xtht1}) and the definition of $G^{\theta,\sigma}_{\alpha_x^{-1}}(x^{t+1},\lambda^{t+1})$, we obtain  that
$$
\begin{array}{rcl}
\|G^{\theta,\sigma}_{\alpha_x^{-1}}(x^{t+1},\lambda^{t+1})\|&=&\alpha_x^{-1}\|(x^{t+1},\lambda^{t+1})-{\rm prox}_{\alpha_x \sigma}((x^{t+1},\lambda^{t+1})-\alpha_x \nabla \theta (x^{t+1},\lambda^{t+1}))\|\\[6pt]
&=&\alpha_x^{-1}\|{\rm prox}_{\alpha_x \sigma}((x^t,\lambda^t)-\alpha_x \xi^t)-{\rm prox}_{\alpha_x \sigma}((x^{t+1},\lambda^{t+1})-\alpha_x \nabla \theta (x^{t+1},\lambda^{t+1}))\|\\[6pt]
&\leq& \alpha_x^{-1}\|(x^{t+1},\lambda^{t+1})-(x^t,\lambda^t)\|+\|\xi^t-\nabla \theta (x^{t+1},\lambda^{t+1})\|\\[6pt]
&\leq& \alpha_x^{-1}\|(x^{t+1},\lambda^{t+1})-(x^t,\lambda^t)\|+\|\xi^t-\nabla \theta (x^{t},\lambda^{t})\|+\|\nabla \theta (x^{t},\lambda^{t})-\nabla \theta (x^{t+1},\lambda^{t+1})\|.
\end{array}
$$
From (\ref{eq:Be0}) and (\ref{eq:xita2}), this implies
$$
\|G^{\theta,\sigma}_{\alpha_x^{-1}}(x^{t+1},\lambda^{t+1})\|\leq (\alpha_x^{-1}+L_{\theta})\|(x^{t+1},\lambda^{t+1})-(x^t,\lambda^t)\|+\rho_0\sqrt{[\|K\|^2+\|B\|^2]}\varepsilon_t.
$$
Therefore, we obtain $\|G^{\theta,\sigma}_{\alpha_x^{-1}}(x^{t+1},\lambda^{t+1})\|\rightarrow 0$ from (\ref{eq:epsilonC}) and (\ref{eq:ccd}). From the definition of $G^{Q,\psi}_{1}(x,y,\lambda)$, one has
$$
\begin{array}{l}
\!\!\|G^{Q,\psi}_{1}(x^{t+1},y^{t+1},\lambda^{t+1})\|\\[6pt]
 \leq \|G^{Q,\psi}_{1}(x^{t+1},y^{t+1},\lambda^{t+1})-G^{Q,\psi}_{1}(x^t,y^{t+1},\lambda^t)\| +\|G^{Q,\psi}_{1}(x^t,y^{t+1},\lambda^t)\|\\[6pt]
= \|G^{Q,\psi}_{1}(x^t,y^{t+1},\lambda^t)\|\\[6pt]
+\|{\rm prox}_{\psi}(y^{t+1}-\nabla h(y^{t+1})+K^Tx^{t+1}+B^T\lambda^{t+1})-{\rm prox}_{\psi}(y^{t+1}-\nabla h(y^{t+1})+K^Tx^{t}+B^T\lambda^{t})\|\\[6pt]
 \leq \varepsilon_t+\|K\|\|x^{t+1}-x^t\|+\|B\|\|\lambda^{t+1}-\lambda^t\|,
\end{array}
$$
which converges to zero when $t \rightarrow \infty$ from (\ref{eq:epsilonC}) and (\ref{eq:ccd}).

Let $(\overline x,\overline y,\overline \lambda)$ be an accumulation point of $\{(x^t,y^t,\lambda^t)\}$. From the Lipschitz continuity of $G^{\theta,\sigma}_{\alpha_x^{-1}}$ and the Lipschitz continuity of $G^{Q,\psi}_{1}$, we obtain
$$
G^{\theta,\sigma}_{\alpha_x^{-1}}(\overline x,\overline \lambda)=0,\quad G^{Q,\psi}_{1}(\overline x,\overline y,\overline \lambda)=0.
$$
From $ G^{Q,\psi}_{1}(\overline x,\overline y,\overline \lambda)=0$, we obtain
\begin{equation}\label{eq:KKTy}
K^T\overline x+B^T\overline \lambda-\nabla h(\overline y) \in \partial \psi (\overline y).
\end{equation}
This implies from the Moreau-Fenchel equality for a proper closed convex function that
$$
\overline y=\nabla \psi_h^*(K^T\overline x+B^T\overline \lambda).
$$
In view of
(\ref{eq:ystar}), $y_*(\overline x,\overline \lambda)=\nabla \psi_h^*(K^T\overline x+B^T\overline \lambda)$, we obtain $
\overline y=y_*(\overline x,\overline \lambda)$.
From (\ref{eq:h02}), we have
$$
\begin{array}{ll}
\nabla \theta (x,\lambda)&\!\!\! =\left
(
\begin{array}{l}
\nabla g(x)+A^T\lambda+Ky_*(x,\lambda)\\
Ax+c+By_*(x,\lambda)
\end{array}
\right).
\end{array}
$$
Thus we have from $G^{\theta,\sigma}_{\alpha_x^{-1}}(\overline x,\overline \lambda)=0$ that
$$
\begin{array}{l}
0\in \nabla g(\overline x)+A^T\lambda+Ky_*(\overline x,\overline \lambda)+\partial \varphi (\overline x),\\[4pt]
0=A\overline x+c+By_*(\overline x,\overline \lambda),
\end{array}
$$
which imply that
\begin{equation}\label{eq:KKTx}
\begin{array}{l}
0\in \nabla g(\overline x)+A^T\lambda+K\overline y+\partial \varphi (\overline x),\\[4pt]
0=A\overline x+c+B\overline y.
\end{array}
\end{equation}
Combining (\ref{eq:KKTy}) and (\ref{eq:KKTx}), we prove (iii). \hfill $\Box$
\begin{remark}\label{remark-about-ac}
The condition $\|y^{t+1}-y_*(t)\|\leq \|y^t-y_*(t)\|$  required by Proposition \ref{prop-4exact} can easily be guaranteed, for example we may use  the gradient ascent method with a constant stepsize.
\end{remark}
\section{A proximal gradient method }
\setcounter{equation}{0}

The following algorithm is based on  getting an approximate stationary point of Problem (\ref{prob-RP}), which is a  proximal gradient multi-step ascent decent method.
\begin{algorithm}\label{alg-30}
Input $(x^0,y^0,\lambda^0) \in \Re^n \times \Re^q \times \Re^m$, $\alpha_x >0$, $\alpha_y(t) >0$ for $t \in \textbf{N}$, positive integers $T>1$ and $N_t>1$ for $t=0,1,2,\ldots,T$\\
\mbox{}\hspace{0.2cm} {\bf for} $t=0,1,2,\ldots,T-1$, {\bf do}\\[4pt]
\mbox{}\hspace{0.5cm} Set $y^{[0]}(t)=y^t$\\
\mbox{}\hspace{0.8cm} {\bf for} $k=0,1,2,\ldots, N_t-1$, {\bf do}\\[4pt]
\mbox{}\hspace{1.2cm} $y^{[k+1]}(t)={\rm prox}_{\alpha_y(t) \psi}\left[y^{[k]}(t)+\alpha_y(t) \left(-\nabla h(y^{[k]}(t)) + K^Tx^t+B^T\lambda^t\right)\right]$\\[4pt]
\mbox{}\hspace{0.8cm} {\bf end for}\\[4pt]
\mbox{}\hspace{0.8cm}Set  $x^{t+1}={\rm prox}_{\alpha_x \varphi}\left[x^t-\alpha_x \left(\nabla g(x^t)+Ky^{t+1}+A^T\lambda^t\right)\right]$\\[4pt]
\mbox{}\hspace{1.2cm} $\lambda^{t+1}=\lambda^t-\alpha_x [Ax^t+By^{t+1}+c]$\\[4pt]
\mbox{}\hspace{0.2cm} {\bf end for}\\
\mbox{}\hspace{0.2cm} Return $(x^{t+1},y^{t+1},\lambda^{t+1})$ for $t=0,1,2,\ldots, T-1$.
\end{algorithm}

In Algorithm \ref{alg-30}, for fixed point $(x^t,\lambda^t)$, a scheme for solving the approximate optimal solution of $y$ by proximal gradient method is given. For simplicity, we denote $y_*(t)=y_*(x^t,\lambda^t)$, since $y \rightarrow f(x,y,\lambda)$ is $\mu$-strongly concave, we have, for any $L>0$, that
$
G^{f,\psi}_L(x^t,y,\lambda^t)=0
$ if and only if $y=y_*(t)$.

We have the following proposition,  whose proof is based on Theorem 10.29 of \cite{Beck2017}.
\begin{proposition}\label{prop-4e0}
Let Assumption  \ref{assum-1e0} be satisfied and
let $\alpha_y(t) \in (0,1/L_h)$. Consider the sequence $\left\{y^{[k]}(t): k=0,1,\ldots,N_t\right\}$ generated by Algorithm \ref{alg-30}, where $t=0,1,\ldots, T$. Then we have for $k=0,\ldots, N_t-1$,
\begin{itemize}
\item[{\rm (a)}]$\left\|y^{[k+1]}(t)-y_*(t)\right\|^2\leq (1-\mu\alpha_y(t))\left\|y^{[k]}(t)-y_*(t)\right\|^2$;
\item[{\rm (b)}]$\left\|y^{[k+1]}(t)-y_*(t)\right\|^2\leq (1-\mu\alpha_y(t))^{k+1}\|y^t-y_*(t)\|^2$;
\item[{\rm (c)}]$\theta_0 (x^t,\lambda^t)-Q\left(x^t,y^{[k+1]}(t),\lambda^t\right)\leq (2\alpha_y(t))^{-1} (1-\mu\alpha_y(t))^{k}\|y^t-y_*(t)\|^2$.
\end{itemize}
\end{proposition}
{\bf Proof}. From the assumptions in this proposition, the function
$$
q^t(y):= -Q(x^t,y,\lambda^t)-\psi(y)=h(y)-(K^Tx^t+B^T\lambda^t)^Ty-\langle \lambda^t, Ax^t+c\rangle
$$
 is $L_h$-smooth and
  $$y \rightarrow -Q(x^t,y,\lambda^t)=q^t(y)+\psi(y)$$
  is
  $\mu$-strongly convex. Noting that for $t=0,1,\ldots, T-1$, the sequence $\{y^{[k]}(t):k=0,1,\ldots,N\}$ satisfies
$$
y^{[k+1]}(t)={\rm prox}_{\alpha_y(t) \psi}\left[y^{[k]}(t)-\alpha_y(t) \nabla q^t (y^{[k]}(t))\right].
$$
 From Lemma \ref{lem2a}, we have for $\alpha_y(t) \in \left(0,1/L_h\right)$,
 \begin{equation}\label{eq10.37e0}
 \begin{array}{rcl}
 Q(x^t,y^{[k+1]}(t),\lambda^t)-\theta_0 (x^t,\lambda^t)& \geq & \displaystyle \frac{1}{2\alpha_y(t)}\|y_*(t)-y^{[k+1]}(t)\|^2-\displaystyle \frac{1}{2\alpha_y(t)}\|y_*(t)-y^{[k]}(t)\|^2\\[16pt]
 & & +\displaystyle \frac{\mu}{2}\|y_*(t)-y^{[k]}(t)\|^2.
 \end{array}
 \end{equation}
Since $\theta_0 (x^t,\lambda^t)=\max_y Q(x^t,y,\lambda^t)$, $ Q(x^t,y^{[k+1]}(t),\lambda^t)-\theta_0 (x^t,\lambda^t)\leq 0$, we have from the above inequality that
$$
\|y_*(t)-y^{[k+1]}(t)\|^2\leq (1-\mu \alpha_y(t)) \|y_*(t)-y^{[k]}(t)\|^2,
$$
which establishes part (a). Part (b) follows from (a) immediately.  To prove part (c),  by (\ref{eq10.37e0}), we have
$$
\begin{array}{rcl}
\theta_0 (x^t,\lambda^t)-Q\left(x^t,y^{[k+1]}(t),\lambda^t\right) &
 \leq & \displaystyle \frac{\alpha_y^{-1}(t)-\mu}{2}\|y_*(t)-y^{[k]}(t)\|^2-\displaystyle \frac{1}{2\alpha_y(t)}\|y_*(t)-y^{[k+1]}(t)\|^2\\[6pt]
 &\leq & \displaystyle \frac{\alpha_y^{-1}(t)-\mu}{2}\|y_*(t)-y^{[k]}(t)\|^2\\[6pt]
& = & (2\alpha_y(t))^{-1} (1-\mu\alpha_y(t))^{k}\|y^t-y_*(t)\|^2,
\end{array}
$$
where part (b) was used in the last inequality. This completes our proof. \hfill $\Box$\\
\begin{theorem}\label{th-este0}
Let Assumption  \ref{assum-1e0} be satisfied and
let $\alpha_y(t) \in \left(0,1/L_h\right)$. Consider the sequence $\left\{y^{[k]}(t): k=0,1,\ldots,N_t\right\}$ generated by Algorithm \ref{alg-30}, and $y^{t+1}=y^{[N_t]}(t)$, where $t=0,1,\ldots, T$. Then one has  for $t=0,1,\ldots,T-1$ that
\begin{equation}\label{eq:innerIneqe0-1}
\begin{array}{l}
\left\|y^{t+1}-y_*(t)\right\|^2\leq \displaystyle (1-\mu\alpha_y(t))^{N_t}\|y^t-y_*(t)\|^2,\\[8pt]
\|\nabla_x Q(x^t,y^{t+1},\lambda^t)-\nabla_x \theta_0 (x^t,\lambda^t)\|^2\leq \displaystyle \|K\|^2 (1-\mu\alpha_y(t))^{N_t}\|y^t-y_*(t)\|^2,\\[8pt]
\|\nabla_{\lambda} Q(x^t,y^{t+1},\lambda^t)-\nabla_{\lambda} \theta_0 (x^t,\lambda^t)\|^2\leq \displaystyle \|B\|^2 (1-\mu\alpha_y(t))^{N_t}\|y^t-y_*(t)\|^2,\\[8pt]
\|G^{f,\psi}_{\alpha_y^{-1}(t)}(x^t,y^{t+1},\lambda^t)\|^2\leq 9\alpha_y^{-2}(t) (1-\mu\alpha_y(t))^{N_t}\|y^t-y_*(t)\|^2
\end{array}
\end{equation}
and
\begin{equation}\label{eq:innerIneqe0}
\begin{array}{l}
\left\|y^{t+1}-y_*(t)\right\|^2\leq \displaystyle \frac{2}{\mu}(1-\mu\alpha_y(t))^{N_t}[\theta_0 (x^t,\lambda^t)-Q(x^t,y^t,\lambda^t)],\\[6pt]
\|\nabla_x Q(x^t,y^{t+1},\lambda^t)-\nabla_x \theta_0 (x^t,\lambda^t)\|^2\leq \displaystyle \frac{2\|K\|^2}{\mu}(1-\mu\alpha_y(t))^{N_t}(\theta_0(x^t,\lambda^t)-Q(x^t,y^t,\lambda^t)),\\[8pt]
\|\nabla_{\lambda} Q(x^t,y^{t+1},\lambda^t)-\nabla_{\lambda} \theta_0 (x^t,\lambda^t)\|^2\leq \displaystyle \frac{2\|B\|^2}{\mu}(1-\mu\alpha_y(t))^{N_t}(\theta_0(x^t,\lambda^t)-Q(x^t,y^t,\lambda^t)),\\[8pt]
\|G^{f,\psi}_{\alpha_y^{-1}(t)}(x^t,y^{t+1},\lambda^t)\|^2\leq \displaystyle \frac{18\alpha_y^{-2}(t)}{\mu}(1-\mu\alpha_y(t))^{N_t}(\theta_0(x^t,\lambda^t)-Q(x^t,y^t,\lambda^t)).
\end{array}
\end{equation}
\end{theorem}
{\bf Proof}. It follows from Proposition \ref{prop-4e0} (b) that
\begin{equation}\label{eq:hh20}
\left\|y^{t+1}(t)-y_*(t)\right\|^2\leq (1-\mu\alpha_y(t))^{N_t}\|y^t-y_*(t)\|^2,
\end{equation}
which is just  the first inequality in (\ref{eq:innerIneqe0-1}).
 From Proposition \ref{prop-1e0}, we have the following equalities
$$
\nabla_x \theta_0 (x^t,\lambda^t)=\nabla_x Q(x^t,y_*(t),\lambda^t)=Ky_*(t)+A^T\lambda^t,\, \nabla_\lambda \theta (x^t,\lambda^t)=\nabla_\lambda Q(x^t,y_*(t),\lambda^t)=Ax^t+By_*(t)+c.
$$
Noting that $\nabla_x Q(x^t,\cdot,\lambda^t)$ and $\nabla_{\lambda} Q(x^t,\cdot,\lambda^t)$ are $\|K\|$-Lipschitz continuous and $\|B\|$-Lipschitz continuous, respectively, we obtain
$$
\|\nabla_x Q(x^t,y^{t+1},\lambda^t)-\nabla_x \theta_0 (x^t,\lambda^t)\|^2=\|\nabla_x Q(x^t,y^{t+1},\lambda^t)-\nabla_x Q(x^t,y_*(t),\lambda^t)\|^2\leq \|K\|^2
\|y^{t+1}-y_*(t)\|^2
$$
and
$$
\|\nabla_{\lambda} Q(x^t,y^{t+1},\lambda^t)-\nabla_{\lambda} \theta (x^t,\lambda^t)\|^2=\|\nabla_z Q(x^t,y^{t+1},\lambda^t)-\nabla_z Q(x^t,y_*(t),\lambda^t)\|^2\leq \|B\|^2
\|y^{t+1}-y_*(t)\|^2.
$$
Thus the second and the third  inequalities in (\ref{eq:innerIneqe0-1}) are obtained  from the first inequality in (\ref{eq:innerIneqe0-1}).

Noting that $G^{f,\psi}_{\alpha_y^{-1}(t)}(x^t,y_*(t),\lambda^t)=0$, we have from Lemma \ref{lem1} that
$$
\|G^{f,\psi}_{\alpha_y^{-1}(t)}(x^t,y^{t+1},\lambda^t)\|^2
=\|G^{f,\psi}_{\alpha_y^{-1}(t)}(x^t,y^{t+1},\lambda^t)-G^{f,\psi}_{\alpha_y^{-1}(t)} (x^t,y_*(t),\lambda^t)\|^2\leq 9\alpha_y^{-2}(t)\|y^{t+1}-y_*(t)\|^2,
$$
and the fourth inequality  in (\ref{eq:innerIneqe0-1}) can be obtained from the first inequality in (\ref{eq:innerIneqe0-1}).

Noting that $y \rightarrow Q(x^t,y,\lambda^t)$ is $\mu$-strongly concave, we have
$$
\theta_0 (x^t,\lambda^t)-Q(x^t,y,\lambda^t) \geq \displaystyle \frac{\mu}{2}\|y-y_*(t)\|^2,
$$
which implies
\begin{equation}\label{eq:hh10}
\|y^t-y_*(t)\|^2 \leq \displaystyle \frac{2}{\mu}[\theta_0 (x^t,\lambda^t)-Q(x^t,y^t,\lambda^t)].
\end{equation}
Thus the inequalities of  (\ref{eq:innerIneqe0}) come   from (\ref{eq:hh10}) and (\ref{eq:hh20}) directly.
The proof is completed.
\hfill $\Box$\\
\begin{proposition}\label{prop-5e0}
Let Assumption  \ref{assum-1e0} be satisfied.
Let $\alpha_y(t) \in \left(0,1/L_h\right)$, $\alpha_x \in \left(0,1/L_{\theta}\right)$  and the sequence $\left\{(x^t,\lambda^t,y^t): t=0,1,\ldots,T\right\}$ be generated by Algorithm \ref{alg-30}. Then for $t=0,\ldots, T-1$,
\begin{equation}\label{eq:keyIneqe0-1}
\begin{array}{rcl}
\varphi(x^{t+1})+\theta (x^{t+1},\lambda^{t+1})& \leq & \theta (x^t,\lambda^t)+\varphi (x^t)-\displaystyle \frac{1-\alpha_x L_{\theta}}{4\alpha_x}\|(x^{t+1},\lambda^{t+1})- (x^t,\lambda^t)\|^2\\[10pt]
& &+\displaystyle \frac{\alpha_x(\|K\|^2+\|B\|^2)}{(1-\alpha_xL_{\theta})}(1-\mu\alpha_y(t))^{N_t}\|y^t-y_*(t)\|^2
\end{array}
\end{equation}
and
\begin{equation}\label{eq:keyIneqe0}
\begin{array}{rcl}
\varphi(x^{t+1})+\theta (x^{t+1},\lambda^{t+1})& \leq & \theta (x^t,\lambda^t)+\varphi (x^t)-\displaystyle \frac{1-\alpha_x L_{\theta}}{4\alpha_x}\|(x^{t+1},\lambda^{t+1})- (x^t,\lambda^t)\|^2\\[10pt]
& & +\displaystyle \frac{2\alpha_x(\|K\|^2+\|B\|^2)}{\mu(1-\alpha_xL_{\theta})}(1-\mu\alpha_y(t))^{N_t}(\theta_0(x^t,\lambda^t)-Q(x^t,y^t,\lambda^t)).
\end{array}
\end{equation}
\end{proposition}
{\bf Proof}. From the definition of $\xi^t$ by (\ref{eq:xit}), one has that $(x^{t+1},\lambda^{t+1})$ can be expressed as
$$
(x^{t+1},\lambda^{t+1})={\rm prox}_{\alpha_x \sigma}((x^t,\lambda^t)-\alpha_x \xi^t),
$$
where $\sigma (x,\lambda)=\varphi (x)$.  Applying Lemma \ref{lem2} with $h=\theta$ and this $\sigma$, we obtain
$$
\begin{array}{rcl}
\varphi (x^{t+1})+\theta (x^{t+1},\lambda^{t+1}) & \leq & \varphi (x^t)+\theta (x^t,\lambda^t)-\displaystyle \frac{1}{2}(\alpha_x^{-1}-L_{\theta})\|(x^{t+1},\lambda^{t+1})- (x^t,\lambda^t)\|^2\\[8pt]
& & +\langle \nabla \theta (x^t,\lambda^t)-\xi^t, (x^{t+1},\lambda^{t+1})- (x^t,\lambda^t)\rangle.
\end{array}
$$
Using the simple inequality $\langle a, b\rangle\leq \|a\|^2/2+\|b\|^2/2$ for vectors $a$ and $b$, we obtain
$$
\begin{array}{rcl}
\varphi(x^{t+1})+\theta (x^{t+1},\lambda^{t+1})&\!\! \leq & \!\! \varphi(x^t)+\theta (x^t,\lambda^t)-\displaystyle \frac{1}{2}(\alpha_x^{-1}-L_{\theta})\|(x^{t+1},\lambda^{t+1})- (x^t,\lambda^t)\|^2\\[6pt]
&\!\! &\!\! +\displaystyle \frac{1}{4}(\alpha_x^{-1}-L_{\theta})\|(x^{t+1},\lambda^{t+1})- (x^t,\lambda^t)\|^2
+\displaystyle \frac{1}{(\alpha_x^{-1}-L_{\theta})}
\|\nabla \theta (x^t,\lambda^t)-\xi^t\|^2.
\end{array}
$$
Then, noting that $\xi^t=\nabla_{x,\lambda}f(x^t,y^{t+1},\lambda^t)$ and
 $$
 \nabla_{x,\lambda} \theta_0 (x^t,\lambda^t)-\nabla_{x,\lambda} Q(x^t,y^{t+1},\lambda^t)=
 \nabla_{x,\lambda} \theta (x^t,\lambda^t)-\xi^t,
 $$
  from (\ref{eq:innerIneqe0-1}) and (\ref{eq:innerIneqe0}) in Theorem \ref{th-este0}, we have
$$
\begin{array}{rcl}
\varphi(x^{t+1})+\theta (x^{t+1},\lambda^{t+1})& \leq & \varphi (x^t)+\theta (x^t,\lambda^t)-\displaystyle \frac{1}{4}(\alpha_x^{-1}-L_{\theta})\|(x^{t+1},\lambda^{t+1})- (x^t,\lambda^t)\|^2\\[8pt]
& &
+\displaystyle \frac{(\|K\|^2+\|B\|^2)}{(\alpha_x^{-1}-L_{\theta})}
(1-\mu\alpha_y(t))^{N_t}\|y^t-y_*(t)\|^2
\end{array}
$$
and
$$
\begin{array}{rcl}
\varphi(x^{t+1})+\theta (x^{t+1},\lambda^{t+1})& \leq & \varphi (x^t)+\theta (x^t,\lambda^t)-\displaystyle \frac{1}{4}(\alpha_x^{-1}-L_{\theta})\|(x^{t+1},\lambda^{t+1})- (x^t,\lambda^t)\|^2\\[8pt]
& & 
+\displaystyle \frac{1}{(\alpha_x^{-1}-L_{\theta})}
\displaystyle \frac{2(\|K\|^2+\|B\|^2)}{\mu}(1-\mu\alpha_y(t))^{N_t}(\theta_0(x^t,\lambda^t)-Q(x^t,y^t,\lambda^t)),
\end{array}
$$
which are just the inequalities (\ref{eq:keyIneqe0-1}) and (\ref{eq:keyIneqe0}). This completes our proof. \hfill $\Box$\\

Let $$\theta_{\varphi} (x,\lambda):=\varphi(x)+\theta (x,\lambda).$$
Then
$$
\begin{array}{rcl}
\theta_{\varphi} (x,\lambda) & = &\varphi(x)+g(x)+\displaystyle \sup_y \left[x^TKy-h(y)-\psi(y)+\langle \lambda,Ax+By+c\rangle\right]\\[6pt]
&= & \varphi_g(x)+\lambda^T(Ax+c)+\psi_h^*(K^Tx+B^T\lambda).
\end{array}
$$
Define
$$
\begin{array}{l}
\chi_0=\displaystyle \frac{1}{\alpha_x(1-L_{\theta}\alpha_x)}, \quad
\chi_1= \displaystyle \frac{(\|K\|^2+\|B\|^2)}{(1-\alpha_x L_{\theta})^2},
\end{array}
$$
and
$$
\delta_t=(1-\mu\alpha_y(t))^{N_t}\|y^t-y_*(t)\|^2,\,\,\,\, \Delta_t=(1-\mu\alpha_y(t))^{N_t}(\theta_0(x^t,\lambda^t)-Q(x^t,y^t,\lambda^t)).
$$
\begin{proposition}\label{prop-main-1e0}
Let Assumptions  \ref{assum-1e0} and  \ref{assum-2e0}  be satisfied.
Let $\alpha_y(t) \in \left(0,1/L_h\right)$, $\alpha_x \in (0,1/L_{\theta})$  and the sequence $\left\{(x^t,\lambda^t,y^t): t=0,1,\ldots,T\right\}$ be generated by Algorithm \ref{alg-30}.
Then there exists an integer $t \in \{0,1,\ldots, T-1\}$ such that
\begin{equation}\label{astatCom0-1}
\begin{array}{l}
\|G^{f,\psi}_{\alpha_y^{-1}(t)}(x^{t+1},y^{t+1},\lambda^{t+1})\|^2\leq 6\chi_0\displaystyle \frac{\mu^2}{\|B\|^2}\displaystyle \frac{\theta (x^0,\lambda^0)-\inf \theta_{\varphi}}{T}+3 \left[\displaystyle \frac{9}{\alpha_y^2(t)}+\displaystyle \frac{2\mu^2 }{\|B\|^2}\chi_1\right]\delta_t,\\[10pt]
 \|G^{f,\varphi}_{\alpha_x^{-1}}(x^{t+1},y^{t+1},\lambda^{t+1})\|^2\leq 4(1+\alpha_xL_g)^2\left\{\chi_0\displaystyle \frac{\theta (x^0,\lambda^0)-\inf \theta_{\varphi}}{T}+\chi_1\delta_t\right\},\\[10pt]
 \|G^{f,0}_{\alpha_x^{-1}}(x^{t+1},y^{t+1},\lambda^{t+1})\|^2\leq 2\chi_0\displaystyle \frac{\theta (x^0,\lambda^0)-\inf \theta_{\varphi}}{T}+
12\chi_1\delta_t
 \end{array}
\end{equation}
and
\begin{equation}\label{astatCom0}
\begin{array}{l}
\|G^{f,\psi}_{\alpha_y^{-1}(t)}(x^{t+1},y^{t+1},\lambda^{t+1})\|^2\leq 6\chi_0\displaystyle \frac{\mu^2}{\|B\|^2}\displaystyle \frac{\theta (x^0,\lambda^0)-\inf \theta_{\varphi}}{T}+6 \left[\displaystyle \frac{9}{\mu\alpha_y^2(t)}+\displaystyle \frac{2\mu }{\|B\|^2}\chi_1\right]\Delta_t,\\[10pt]
 \|G^{f,\varphi}_{\alpha_x^{-1}}(x^{t+1},y^{t+1},\lambda^{t+1})\|^2\leq 4(1+\alpha_xL_g)^2\left\{\chi_0\displaystyle \frac{\theta (x^0,\lambda^0)-\inf \theta_{\varphi}}{T}+\frac{2}{\mu}\chi_1\Delta_t\right\},\\[10pt]
 \|G^{f,0}_{\alpha_x^{-1}}(x^{t+1},y^{t+1},\lambda^{t+1})\|^2\leq 12 \left[\chi_0\displaystyle \frac{\theta (x^0,\lambda^0)-\inf \theta_{\varphi}}{T}+
\displaystyle \frac{2}{\mu}\chi_1\Delta_t\right].
 \end{array}
\end{equation}
\end{proposition}
{\bf Proof}. We only prove (\ref{astatCom0}). The relation (\ref{astatCom0-1}) can be proved similarly. In view of Assumptions  \ref{assum-1e0} and \ref{assum-2e0}, summing (\ref{eq:keyIneqe0})
 over $t=0,1,\ldots, T-1$,
  we obtain
 $$
\begin{array}{rcl}
\displaystyle \sum_{t=0}^{T-1}\|(x^{t+1},\lambda^{t+1})-(x^{t},\lambda^{t})\|^2
 &\leq& \displaystyle \frac{4\alpha_x}{1-\alpha_xL_{\theta}}[\theta_{\varphi} (x^0,\lambda^0)-\theta_{\varphi} (x^T,\lambda^T)] +\displaystyle \frac{8\alpha_x^2(\|K\|^2+\|B\|^2)}{\mu (1-\alpha_xL_{\theta})^2}\sum_{t=0}^{T-1}\Delta_t\\[8pt]
  &\leq& \displaystyle \frac{4\alpha_x}{1-\alpha_xL_{\theta}}[\theta_{\varphi} (x^0,\lambda^0)-\inf \theta_\varphi]+\displaystyle \frac{8\alpha_x^2(\|K\|^2+\|B\|^2)}{\mu (1-\alpha_xL_{\theta})^2}\sum_{t=0}^{T-1}\Delta_t.
\end{array}
$$
Thus there exists an integer $t \in \{0,1,\ldots, T-1\}$ such that
\begin{equation}\label{eq:xd0}
\begin{array}{ll}
\|(x^{t+1},\lambda^{t+1})-(x^{t},\lambda^{t})\|^2 \leq \displaystyle \frac{4\alpha_x}{1-\alpha_xL_{\theta}}\displaystyle \frac{\theta (x^0,\lambda^0)-\inf \theta_{\varphi}}{T} +
\displaystyle \frac{8\alpha_x^2(\|K\|^2+\|B\|^2)}{\mu (1-\alpha_xL_{\theta})^2}\Delta_t.
\end{array}
\end{equation}
Since
$$
y_*(t)={\rm argmax}\left\{(x^{t})^TKy-h(y)-\psi (y)+\langle \lambda^t, Ax^t+By+c\rangle\right\},
$$
we have
$$
0 \in -[ K^Tx^t+B^T\lambda^t]+\nabla h(y_*(t))+\partial \psi (y_*(t)),
$$
implying that
$$
y_*(t)={\rm prox}_{\alpha_y(t)\psi}(y_*(t)+\alpha_y(t)( -\nabla h(y_*(t))+K^Tx^t+B^T\lambda^t))=
{\rm prox}_{\alpha_y(t)\psi}(y_*(t)+\alpha_y(t)\nabla_y (x^t,\lambda^t,y_*(t))).
$$
Then we get
$$
\begin{array}{ll}
&\!\!\!\!\|G^{f,\psi}_{\alpha_y^{-1}(t)}(x^{t+1},y^{t+1},\lambda^{t+1})\|\\[6pt]
&=\|\alpha_y^{-1}(t)[y^{t+1}-{\rm prox}_{\alpha_y(t)\psi}(y^{t+1}+\alpha_y(t) \nabla_y f(x^{t+1},\lambda^{t+1},y^{t+1}))]\|\\[6pt]
&= \|\alpha_y^{-1}(t)[y^{t+1}-y_*(t)]- \alpha_y^{-1}(t){\rm prox}_{\alpha_y(t)\psi}(y^{t+1}+\alpha_y\nabla_y f(x^{t+1},\lambda^{t+1},y^{t+1}))\\[6pt]
&\quad  +\alpha_y^{-1}(t){\rm prox}_{\alpha_y(t)\psi}(y_*(t)+\alpha_y(t)\nabla_y (x^t,\lambda^t,y_*(t)))\|\\[6pt]
& \leq 2\alpha_y^{-1}(t) \|y^{t+1}-y_*(t)\|+\|\nabla_y f(x^{t+1},\lambda^{t+1},y^{t+1})-\nabla_y f(x^t,\lambda^t,y_*(t))\|\\[6pt]
&=2\alpha_y^{-1}(t) \|y^{t+1}-y_*(t)\|
+\|K^T(x^{t+1}-x^t)+B^T(\lambda^{t+1}-\lambda^t)-[\nabla h(y^{t+1})-\nabla h(y_*(t))]\|\\[6pt]
& \leq (L_h+2\alpha_y^{-1}(t)) \|y^{t+1}-y_*(t)\|+\|B\|\|\lambda^{t+1}-\lambda^t\|+\|K\| \|x^{t+1}-x^t\|,
\end{array}
$$
which implies
$$
\|G^{f,\psi}_{\alpha_y^{-1}(t)}(x^{t+1},y^{t+1},\lambda^{t+1})\|^2\leq 3(L_h+2\alpha_y^{-1}(t))^2 \|y^{t+1}-y_*(t)\|^2+3 \max[\|B\|^2,\|K\|^2]\|(x^{t+1},\lambda^{t+1})-(x^{t},\lambda^{t})\|^2.
$$
From Theorem \ref{th-este0} and the relation (\ref{eq:xd0}), we obtain
$$
\begin{array}{l}
\!\!\|G^{f,\psi}_{\alpha_y^{-1}(t)}(x^{t+1},y^{t+1},\lambda^{t+1})\|^2\\[8pt] \leq \displaystyle \frac{6[L_h+2\alpha_y^{-1}(t)]^2 }{\mu}\Delta_t+3 \max[\|B\|^2,\|K\|^2]\left\{ \displaystyle \frac{4\alpha_x}{1-\alpha_xL_{\theta}}\displaystyle \frac{\theta (x^0,\lambda^0)-\inf \theta_{\varphi}}{T}+
\displaystyle \frac{8\alpha_x^2(\|K\|^2+\|B\|^2)}{\mu (1-\alpha_xL_{\theta})^2}\Delta_t\right \}.
\end{array}
$$
From the definition of $L_{\theta}$ and the choice of
  $\alpha_x \in (0,1/L_{\theta})$, we obtain
  $$
  \max[\|B\|^2,\|K\|^2]\alpha_x^2<[ \|B\|^2+\|K\|^2]\alpha_x^2<[ \|B\|^2+\|K\|^2]\displaystyle \frac{\mu^2}{\gamma^2}
  \leq \displaystyle \frac{\mu^2}{2\|B\|^2}.
  $$
  Thus we obtain from $\alpha_y(t) \in (0, 1/L_h)$ that $L_h+2\alpha_y^{-1}(t)< 3\alpha_y^{-1}(t)$ and
\begin{equation}\label{eq:es10}
\begin{array}{l}
\!\!\!\!\|G^{f,\psi}_{\alpha_y^{-1}(t)}(x^{t+1},y^{t+1},\lambda^{t+1})\|^2 \\[8pt]
 \leq \displaystyle \frac{6[L_h+2\alpha_y^{-1}(t)]^2 }{\mu}\Delta_t
 +\displaystyle \frac{\mu^2}{\|B\|^2}\left\{ \displaystyle \frac{6}{\alpha_x(1-\alpha_xL_{\theta})}\displaystyle \frac{\theta (x^0,\lambda^0)-\inf \theta_{\varphi}}{T}+
\displaystyle \frac{12(\|K\|^2+\|B\|^2)}{\mu (1-\alpha_xL_{\theta})^2}\Delta_t\right \}\\[16pt]
\leq 6 \left[\displaystyle \frac{9}{\mu\alpha_y^2(t)}+\displaystyle \frac{2\mu }{\|B\|^2}\chi_1\right]\Delta_t+6\chi_0\displaystyle \frac{\mu^2}{\|B\|^2}\displaystyle \frac{\theta (x^0,\lambda^0)-\inf \theta_{\varphi}}{T}.
\end{array}
\end{equation}
From Algorithm \ref{alg-30}, we can express $(x^{t+1},\lambda^{t+1})$ as
$$
(x^{t+1},\lambda^{t+1})={\rm prox}_{\alpha_x \sigma}((x^t,\lambda^t)-\alpha_x \xi^t),
$$
with $\xi^t$ defined by (\ref{eq:xit}).  Thus we have
$$
\begin{array}{l}
\!\!\!\!\|G^{f,\varphi}_{\alpha_x^{-1}}(x^{t+1},y^{t+1},\lambda^{t+1})\| \\[6pt]
=\|\alpha_x^{-1}[x^{t+1}-{\rm prox}_{\alpha_x \varphi}(x^{t+1}-\alpha_x \nabla_x f(x^{t+1},\lambda^{t+1},y^{t+1}))\|\\[6pt]
=\alpha_x^{-1}\|{\rm prox}_{\alpha_x \varphi}(x^t-\alpha_x \nabla_x f(x^t,y^{t+1},\lambda^{t}))-{\rm prox}_{\alpha_x \varphi}(x^{t+1}-\alpha_x \nabla_x f(x^{t+1},y^{t+1},\lambda^{t+1}))\|\\[6pt]
=\alpha_x^{-1}\|{\rm prox}_{\alpha_x \varphi}(x^t-\alpha_x[\nabla g(x^t)+Ky^{t+1}+A^T\lambda^{t}])
\\[6pt]\quad -{\rm prox}_{\alpha_x \varphi}(x^{t+1}-\alpha_x [\nabla g(x^{t+1})+Ky^{t+1}+A^T\lambda^{t+1}])\|\\[6pt]
 \leq \alpha_x^{-1}\|x^{t+1}-x^t\|+L_g \|x^{t+1}-x^t\|+\|A\|\|\lambda^{t+1}-\lambda^t\|.
\end{array}
$$
Then we obtain from $\alpha_x^{-2}>2 \|A\|^2$, $\max\{(\alpha_x^{-1}+L_g)^2,\|A\|^2\} = (\alpha_x^{-1}+L_g)^2$, that
\begin{equation}\label{eq:es2}
\begin{array}{l}
\!\!\!\!\|G^{f,\varphi}_{\alpha_x^{-1}}(x^{t+1},y^{t+1},\lambda^{t+1})\|^2\\
 \leq 2(\alpha_x^{-1}+L_g)^2\|x^{t+1}-x^t\|^2+2\|A\|^2\|\lambda^{t+1}-\lambda^t\|^2\\[6pt]
 \leq \max\{(\alpha_x^{-1}+L_g)^2,\|A\|^2\} \left\{ \displaystyle \frac{4\alpha_x}{1-\alpha_xL_{\theta}}\displaystyle \frac{\theta (x^0,\lambda^0)-\inf \theta_{\varphi}}{T}+
\displaystyle \frac{8\alpha_x^2(\|K\|^2+\|B\|^2)}{\mu (1-\alpha_xL_{\theta})^2}\Delta_t\right \}\\[16pt]
\leq (1+\alpha_xL_g)^2\left\{ \displaystyle \frac{4}{\alpha_x(1-\alpha_xL_{\theta})}\displaystyle \frac{\theta (x^0,\lambda^0)-\inf \theta_{\varphi}}{T}+
\displaystyle \frac{8(\|K\|^2+\|B\|^2)}{\mu (1-\alpha_xL_{\theta})^2}\Delta_t\right \}\\[16pt]
=4(1+\alpha_xL_g)^2\left\{\chi_0\displaystyle \frac{\theta (x^0,\lambda^0)-\inf \theta_{\varphi}}{T}+\chi_1\displaystyle \frac{2}{\mu}\Delta_t\right\}.
\end{array}
\end{equation}
From the updating formula
$$
{Ax^t+By^{t+1}+c=\displaystyle \frac{\lambda^{t+1}-\lambda^t}{\alpha_x},}
$$
we have that
$$
G^{f,0}_{\alpha_x^{-1}}(x^{t+1},y^{t+1},\lambda^{t+1})=Ax^{t+1}+By^{t+1}+c=A(x^{t+1}-x^t)+\displaystyle {\frac{\lambda^{t+1}-\lambda^t}{\alpha_x}}.
$$
From the relation $\alpha_x^{-2}>2 \|A\|^2$, we obtain from (\ref{eq:xd0}) that
$$
\begin{array}{rcl}
\|G^{f,0}_{\alpha_x^{-1}}(x^{t+1},y^{t+1},\lambda^{t+1})\|^2
& = &\|Ax^{t+1}+By^{t+1}+c\|^2 \\[6pt]
& \leq & [\|A\|\|x^{t+1}-x^t\|+\alpha_x^{-1}\|\lambda^t-\lambda^{t+1}\|]^2\\[6pt]
&\leq & 2[\|A\|^2\|x^{t+1}-x^t\|^2+\alpha_x^{-2}\|\lambda^t-\lambda^{t+1}\|^2]\\[6pt]
&\leq & 3\alpha_x^{-2}\|(x^{t+1},\lambda^{t+1})-(x^{t},\lambda^{t})\|^2\\[6pt]
& \leq & \displaystyle \frac{12}{\alpha_x(1-\alpha_xL_{\theta})}\displaystyle \frac{\theta (x^0,\lambda^0)-\inf \theta_{\varphi}}{T}+
\displaystyle \frac{24(\|K\|^2+\|B\|^2)}{\mu (1-\alpha_xL_{\theta})^2}\Delta_t.\\[10pt]
& = & 12 \chi_0\displaystyle \frac{\theta (x^0,\lambda^0)-\inf \theta_{\varphi}}{T}+
\displaystyle \frac{24}{\mu}\chi_1\Delta_t.
\end{array}
$$
The proof is completed. \hfill $\Box$\\

Noting that Proposition \ref{prop-main-1e0} provides the estimates for $\|G^{f,0}_{\alpha_x^{-1}}(x^{t+1},y^{t+1},\lambda^{t+1})\|$, {\it i.e.}, the estimates for $\|A x^{t+1}+B y^{t+1}+c\|$, we may modify $(x^{t+1},y^{t+1})$ to satisfy the equality constraint $Ax+By+c=0$. This needs the following assumption, which is not hard to satisfy.
\begin{assumption}\label{assum-AB}
Suppose that $A$ and $B$ satisfy that $[A\,\, B]$ is of full row rank.
\end{assumption}

It follows from Assumption \ref{assum-AB} that $AA^T+BB^T$ is positively definite.
We give the expression of the projection of  a point of $\Re^n \times \Re^m$ onto $C$ in the following lemma.
\begin{lemma}\label{lem:prj}
Let Assumption \ref{assum-AB} be satisfied. Then for any $(\widetilde x,\widetilde y)$, its projection on $C$, denoted by $\Pi_C(\widetilde x,\widetilde y)$, is given by
\begin{equation}\label{eq:pjexp}
\Pi_C(\widetilde x,\widetilde y)=(\widetilde x,\widetilde y)-(A^T\widetilde \zeta,B^T\widetilde \zeta)
\end{equation}
with
$$
\widetilde \zeta=(AA^T+BB^T)^{-1}(A\widetilde x+B\widetilde y+c).
$$
\end{lemma}
Define
\begin{equation}\label{eq:kp1}
(\widetilde x^{t+1},\widetilde y^{t+1})=\Pi_C(x^{t+1},y^{t+1}).
\end{equation}
Let
$$
\omega_y(t)=[L_g+2\alpha_y(t)^{-1}]\|A^T[AA^T+BB^T]^{-1}\|+\|K\|\|B^T[AA^T+BB^T]^{-1}\|
$$
and
$$
\omega_x=[L_g+2\alpha_x^{-1}]\|A^T[AA^T+BB^T]^{-1}\|+\|K\|\|B^T[AA^T+BB^T]^{-1}\|.
$$
Then we have the following conclusion.
\begin{proposition}\label{prop-prj0}
Let Assumptions  \ref{assum-1e0},  \ref{assum-2e0} and \ref{assum-AB} be satisfied.
Let $\alpha_y(t) \in \left(0,1/L_h\right)$, $\alpha_x \in (0,1/L_{\theta})$  and the sequence $\left\{(x^t,\lambda^t,y^t): t=0,1,\ldots,T\right\}$ be generated by Algorithm \ref{alg-30}.
Then there exists an integer $t \in \{0,1,\ldots, T-1\}$ such that
$$
 \|G^{f,0}_{\alpha_x^{-1}}(\widetilde x^{t+1},\widetilde y^{t+1},\lambda^{t+1})\|^2 =  \|A\widetilde x^{t+1}+B\widetilde y^{t+1}+c\|^2=0,
$$
\begin{equation}\label{onprj0-1}
\begin{array}{rcl}
\|G^{f,\psi}_{\alpha_y^{-1}(t)}(\widetilde x^{t+1},\widetilde y^{t+1},\lambda^{t+1})\|^2 & \leq &  12 \chi_0 \left[\displaystyle \frac{\mu^2}{\|B\|^2}+2 \omega_y(t)^2\right]\displaystyle \frac{\theta (x^0,\lambda^0)-\inf \theta_{\varphi}}{T}\\[12pt]
& \quad & +6\left[\displaystyle \frac{9}{\alpha_y^2(t)}+\displaystyle \frac{2\mu^2 }{\|B\|^2}\chi_1+\displaystyle
4\chi_1\omega_y(t)^2\right]\delta_t,\\[16pt]
 \|G^{f,\varphi}_{\alpha_x^{-1}}(\widetilde x^{t+1},\widetilde y^{t+1},\lambda^{t+1})\|^2 & \leq & 8\chi_0\left[(1+\alpha_xL_g)^2+3 \omega_x^2\right]
  \displaystyle \frac{\theta (x^0,\lambda^0)-\inf \theta_{\varphi}}{T}\\[12pt]
  & \quad & +8\chi_1\left[(1+\alpha_xL_g)^2 +\displaystyle
3\omega_x^2\right]\delta_t
 \end{array}
\end{equation}
and
\begin{equation}\label{onprj0}
\begin{array}{rcl}
\|G^{f,\psi}_{\alpha_y^{-1}(t)}(\widetilde x^{t+1},\widetilde y^{t+1},\lambda^{t+1})\|^2 & \leq & 12 \chi_0 \left[\displaystyle \frac{\mu^2}{\|B\|^2}+2 \omega_y(t)^2\right]\displaystyle \frac{\theta (x^0,\lambda^0)-\inf \theta_{\varphi}}{T}\\[12pt]
& \quad & +12 \left[\displaystyle \frac{9}{\mu\alpha_y^2(t)}+\displaystyle \frac{2\mu }{\|B\|^2}\chi_1+\displaystyle
\frac{4}{\mu}\chi_1\omega_y(t)^2\right]\Delta_t,\\[16pt]
 \|G^{f,\varphi}_{\alpha_x^{-1}}(\widetilde x^{t+1},\widetilde y^{t+1},\lambda^{t+1})\|^2 & \leq &  8\chi_0\left[(1+\alpha_xL_g)^2+3 \omega_x^2\right]
  \displaystyle \frac{\theta (x^0,\lambda^0)-\inf \theta_{\varphi}}{T}\\[12pt]
  & \quad & + \displaystyle \frac{16}{\mu}\chi_1\left[(1+\alpha_xL_g)^2 +\displaystyle
3\omega_x^2\right]\Delta_t.
 \end{array}
\end{equation}
\end{proposition}
{\bf Proof}. The equality
$\|A\widetilde x^{t+1}+B\widetilde y^{t+1}+c\|^2=0$ comes from the definition (\ref{eq:kp1}). Here we only prove
(\ref{onprj0}) and (\ref{onprj0-1}) can be proved similarly. Let
$$
\widetilde \zeta^{t+1}=(AA^T+BB^T)^{-1}(A x^{t+1}+B y^{t+1}+c).
$$
Then, from Lemma \ref{lem:prj}, we obtain
$$
\begin{array}{l}
\!\!\!\!\|G^{f,\psi}_{\alpha_y^{-1}(t)}(\widetilde x^{t+1},\widetilde y^{t+1},\lambda^{t+1})\|\\[12pt]
 \leq
\|G^{f,\psi}_{\alpha_y^{-1}(t)}(x^{t+1}, y^{t+1},\lambda^{t+1})\|+\|G^{f,\psi}_{\alpha_y^{-1}(t)}(\widetilde x^{t+1},\widetilde y^{t+1},\lambda^{t+1})-G^{f,\psi}_{\alpha_y^{-1}(t)}(x^{t+1}, y^{t+1},\lambda^{t+1})\|\\[12pt]
\leq \|G^{f,\psi}_{\alpha_y^{-1}(t)}(x^{t+1}, y^{t+1},\lambda^{t+1})\|+2\alpha_y(t)^{-1}\|\widetilde x^{t+1}-x^{t+1}\|+L_g
\|\widetilde x^{t+1}-x^{t+1}\|+\|K\|\|\widetilde y^{t+1}-y^{t+1}\|\\[12pt]
=\|G^{f,\psi}_{\alpha_y^{-1}(t)}(x^{t+1}, y^{t+1},\lambda^{t+1})\|+(L_g+2\alpha_y(t)^{-1})\|A^T\widetilde \zeta^{t+1}\|+\|K\|\|B^T\widetilde \zeta^{t+1}\|\\[12pt]
\leq \|G^{f,\psi}_{\alpha_y^{-1}(t)}(x^{t+1}, y^{t+1},\lambda^{t+1})\|+\omega_y(t)\|Ax^{t+1}+By^{t+1}+c\|.
\end{array}
$$
Thus we have from Proposition \ref{prop-main-1e0} that
$$
\begin{array}{rcl}
\|G^{f,\psi}_{\alpha_y^{-1}(t)}(\widetilde x^{t+1},\widetilde y^{t+1},\lambda^{t+1})\|^2&\!\!
 \leq &\!\! 2\|G^{f,\psi}_{\alpha_y^{-1}(t)}(x^{t+1}, y^{t+1},\lambda^{t+1})\|^2+2\omega_y(t)^2\|Ax^{t+1}+By^{t+1}+c\|^2\\[12pt]
 &\!\!\leq &\!\! 12\chi_0\displaystyle \frac{\mu^2}{\|B\|^2}\displaystyle \frac{\theta (x^0,\lambda^0)-\inf \theta_{\varphi}}{T}+12 \left[\displaystyle \frac{9}{\mu\alpha_y^2(t)}+\displaystyle \frac{2\mu }{\|B\|^2}\chi_1\right]\Delta_t\\[10pt]
 &\!\! &\!\! +2\omega_y(t)^2  \times \left\{ \displaystyle \frac{12}{\alpha_x(1-\alpha_xL_{\theta})}\displaystyle \frac{\theta (x^0,\lambda^0)-\inf \theta_{\varphi}}{T}+
\displaystyle \frac{24(\|K\|^2+\|B\|^2)}{\mu (1-\alpha_xL_{\theta})^2}\Delta_t \right\},
 \end{array}
$$
which yields the first inequality in (\ref{onprj0}).  The second inequality in (\ref{onprj0}) can be proved in the same way.
\hfill $\Box$

For developing the iteration complexity of Algorithm \ref{alg-30}, we need the following extra assumptions.
\begin{assumption}\label{assum-domain-B}
Suppose that there exists a constant $\beta_1>0$ such that
$$
\mbox{dom }\psi \subset \beta_1\textbf{B},
$$
where $\mbox{dom }\psi=\{y: \psi(y) <+\infty\}$ is the effective domain of $\psi$.
\end{assumption}
\begin{remark}\label{Y-constraint}
If $Y \subset \Re^m$ is an nonempty convex compact set with $Y \subset \beta_1\textbf{B}$ for some $\beta_1>0$, then
$$\psi (y)=\delta_Y(y)=\left\{
\begin{array}{ll}
0 & y \in Y,\\[4pt]
+\infty & y \notin Y
\end{array}
\right.
$$
satisfies Assumption
\ref{assum-domain-B}.
\end{remark}
\begin{assumption}\label{assum-3e0}
Suppose that there exists a constant $\omega_1>0$ such that
$$
 \theta_0 (x^t,\lambda^t)-Q(x^t,y^t,\lambda^t) \leq \omega_1
$$
for $t=0,1,\ldots,T$.
\end{assumption}
Define
\begin{equation}\label{eq:wy}
\omega_y=[L_g+2\alpha_y^{-1}]\|A^T[AA^T+BB^T]^{-1}\|+\|K\|\|B^T[AA^T+BB^T]^{-1}\|
\end{equation}
and
\begin{equation}\label{eq:n12}
\begin{array}{l}
\gamma_1=\max\left\{ 6\left[\displaystyle \frac{9}{\alpha_y^2}+\displaystyle \frac{2\mu^2 }{\|B\|^2}\chi_1+\displaystyle
4\chi_1\omega_y^2\right], 8\chi_1\left[(1+\alpha_xL_g)^2 +\displaystyle
3\omega_x^2\right]\right\},\\[18pt]
\gamma_2=\max\left\{ 12 \chi_0 \left[\displaystyle \frac{\mu^2}{\|B\|^2}+2 \omega_y^2\right], 8\chi_0\left[(1+\alpha_xL_g)^2+3 \omega_x^2\right]\right\}.
\end{array}
\end{equation}
\begin{theorem}\label{th-main-1e0}
Let Assumptions  \ref{assum-1e0}, \ref{assum-2e0} and \ref{assum-AB} be satisfied. Let $N_t\equiv N$  be a constant positive integer for $t\in \textbf{N}$.
Let $\alpha_y(t) \equiv \alpha_y\in (0,1/L_h)$, $\alpha_x \in (0,1/L_{\theta})$  and the sequence $\left\{(x^t,\lambda^t,y^t): t=0,1,\ldots,T\right\}$ be generated by Algorithm \ref{alg-30}.
For any $\epsilon>0$, we have the following conclusions:
\begin{itemize}
\item[(i)]If  Assumption  \ref{assum-domain-B} is satisfied and
$$
N\geq  \displaystyle \frac{1}{-\log (1-\mu\alpha_y)}\left[\log(8\gamma_1\beta_1^2)+2\log \displaystyle \frac{1}{\epsilon}\right], \quad T \geq \displaystyle \frac{2\gamma_2(\theta_{\varphi} (x^0,\lambda^0)-\inf \theta_{\varphi} )}{\epsilon^2},
$$
then there exists an integer $t \in \{0,1,\ldots, T-1\}$ such that $(\widetilde x^{t+1},\widetilde y^{t+1},\lambda^{t+1})$ satisfies
$$
\|G^{f,\psi}_{\alpha_y^{-1}}(\widetilde x^{t+1},\widetilde y^{t+1},\lambda^{t+1})\|
\leq \epsilon,\,\,
\|G^{f,\varphi}_{\alpha_x^{-1}}(\widetilde x^{t+1},\widetilde y^{t+1},\lambda^{t+1})\|\leq \epsilon,\,
A\widetilde x^{t+1}+B\widetilde y^{t+1}+c=0.
$$
\item[(ii)]If  Assumption  \ref{assum-3e0} is satisfied and
$$
N\geq  \displaystyle \frac{1}{-\log (1-\mu\alpha_y)}\left[\log\left(\displaystyle \frac{4\gamma_1\omega_1}{\mu}\right)+2\log \displaystyle \frac{1}{\epsilon}\right], \quad  T \geq \displaystyle \frac{2\gamma_2(\theta_{\varphi} (x^0,\lambda^0)-\inf \theta_{\varphi} )}{\epsilon^2},
$$
then there exists an integer $t \in \{0,1,\ldots, T-1\}$ such that
$$
\|G^{f,\psi}_{\alpha_y^{-1}}(\widetilde x^{t+1},\widetilde y^{t+1},\lambda^{t+1})\|
\leq \epsilon,\,\,
\|G^{f,\varphi}_{\alpha_x^{-1}}(\widetilde x^{t+1},\widetilde y^{t+1},\lambda^{t+1})\|\leq \epsilon,\,
A\widetilde x^{t+1}+B\widetilde y^{t+1}+c=0.
$$
\end{itemize}
\end{theorem}
{\bf Proof}. When $\alpha_y(t)\equiv \alpha_y\in
(0,1/(L_h+\beta_0))$ is chosen, it is a constant independent of $t$. In this case, $\omega_y(t)\equiv \omega_y$, where $\omega_y$ is defined by (\ref{eq:wy}). Thus, using $\gamma_1$ and $\gamma_2$ defined by (\ref{eq:n12}), $N_t\equiv N \in \textbf{N}$, $\|y^t-y_*(t)\|\leq 2 \beta_1$ by  Assumption \ref{assum-domain-B},
 we have from (\ref{onprj0-1}) of
Proposition \ref{prop-prj0} that there  exists an integer $t \in \{0,1,\ldots, T-1\}$,
$$
 \|G^{f,0}_{\alpha_x^{-1}}(\widetilde x^{t+1},\widetilde y^{t+1},\lambda^{t+1})\|^2 =  \|A\widetilde x^{t+1}+B\widetilde y^{t+1}+c\|^2=0,
$$
\begin{equation}\label{onprj0-1h}
\begin{array}{rcl}
\|G^{f,\psi}_{\alpha_y^{-1}}(\widetilde x^{t+1},\widetilde y^{t+1},\lambda^{t+1})\|^2 & \leq & \gamma_2\displaystyle \frac{\theta (x^0,\lambda^0)-\inf \theta_{\varphi}}{T}+(1-\mu\alpha_y)^{N}4\beta_1^2,\\[8pt]
 \|G^{f,\varphi}_{\alpha_x^{-1}}(\widetilde x^{t+1},\widetilde y^{t+1},\lambda^{t+1})\|^2 & \leq  & \gamma_2
  \displaystyle \frac{\theta (x^0,\lambda^0)-\inf \theta_{\varphi}}{T}+\gamma_1(1-\mu\alpha_y)^{N}4\beta_1^2.
 \end{array}
\end{equation}
Instead of Assumption \ref{assum-domain-B}, if Assumption \ref{assum-3e0} holds, then $
 \theta_0 (x^t,\lambda^t)-Q(x^t,y^t,\lambda^t) \leq \omega_1$, and  we have from (\ref{onprj0-1}) of
Proposition \ref{prop-prj0} that there  exists an integer $t \in \{0,1,\ldots, T-1\}$,
\begin{equation}\label{onprj0h}
\begin{array}{rcl}
\|G^{f,\psi}_{\alpha_y^{-1}}(\widetilde x^{t+1},\widetilde y^{t+1},\lambda^{t+1})\|^2 & \leq &  \gamma_2\displaystyle \frac{\theta (x^0,\lambda^0)-\inf \theta_{\varphi}}{T}+\displaystyle \frac{2\gamma_1\omega_1}{\mu}(1-\mu\alpha_y)^{N};\\[8pt]
 \|G^{f,\varphi}_{\alpha_x^{-1}}(\widetilde x^{t+1},\widetilde y^{t+1},\lambda^{t+1})\|^2 & \leq  & \gamma_2\displaystyle \frac{\theta (x^0,\lambda^0)-\inf \theta_{\varphi}}{T}+\displaystyle\frac{2\gamma_1\omega_1}{\mu}(1-\mu\alpha_y)^{N}.
 \end{array}
\end{equation}
If  Assumption  \ref{assum-domain-B} is satisfied and
$$
N\geq  \displaystyle \frac{1}{-\log (1-\mu\alpha_y)}\left[\log(8\gamma_1\beta_1^2)+2\log \displaystyle \frac{1}{\epsilon}\right], \quad T \geq \displaystyle \frac{2\gamma_2(\theta_{\varphi} (x^0,\lambda^0)-\inf \theta_{\varphi} )}{\epsilon^2},
$$
we have
$$
\gamma_2\displaystyle \frac{\theta (x^0,\lambda^0)-\inf \theta_{\varphi}}{T}\leq
\displaystyle \frac{\epsilon^2}{2}, \,\,\,\,\, (1-\mu\alpha_y)^{N}4\beta_1^2\leq
\displaystyle \frac{\epsilon^2}{2}.
$$
It follows from this and (\ref{onprj0-1h}) that
$$
\|G^{f,\psi}_{\alpha_y^{-1}}(\widetilde x^{t+1},\widetilde y^{t+1},\lambda^{t+1})\|^2\leq \varepsilon^2,\ \
\|G^{f,\varphi}_{\alpha_x^{-1}}(\widetilde x^{t+1},\widetilde y^{t+1},\lambda^{t+1})\|^2\leq \varepsilon^2,
$$
which implies the truth of (i).   If  Assumption  \ref{assum-3e0} is satisfied and
$$
N\geq  \displaystyle \frac{1}{-\log (1-\mu\alpha_y)}\left[\log\left(\displaystyle \frac{4\gamma_1\omega_1}{\mu}\right)+2\log \displaystyle \frac{1}{\epsilon}\right], \quad  T \geq \displaystyle \frac{2\gamma_2(\theta_{\varphi} (x^0,\lambda^0)-\inf \theta_{\varphi} )}{\epsilon^2},
$$
we have
$$
\gamma_2\displaystyle \frac{\theta (x^0,\lambda^0)-\inf \theta_{\varphi}}{T}\leq
\displaystyle \frac{\epsilon^2}{2}, \,\,\,\,\, \displaystyle \frac{2\gamma_1\omega_1}{\mu}(1-\mu\alpha_y)^{N}\leq
\displaystyle \frac{\epsilon^2}{2},
$$
which with (\ref{onprj0h}) proves (ii). The proof is completed. \hfill $\Box$

Theorem \ref{th-main-1e0} tells us, for the structured linearly constrained  nonsmooth minimax problem (\ref{minimax}), that the multi-step
ascent descent proximal gradient method Algorithm \ref{alg-30} can find an $\epsilon$-stationary point  in ${\cal O}\left(\epsilon^{-2}\log  \epsilon^{-1}  \right)$ iterations if $N$ and $T$ are chosen as in Theorem \ref{th-main-1e0}.
\section{Applications}
In this section, we  apply the proposed algorithm to solve three classes of problems. The first one is generalized absolute value equations (GAVE), which can be converted to a smooth convex-concave minimax problem; i.e., $\varphi\equiv0$ and $\psi\equiv0$. We discuss the numerical algorithm for solving GAVE and test the algorithm for three examples. Secondly, we consider a linear regression problem,  which is a smooth strongly-convex-strongly-concave minimax problem, and we perform the effect of Algorithm \ref{alg-30} with different dimensions and constraints. Thirdly, Algorithm \ref{alg-30} is used to solve generalized linear projection equations, which is equivalent to a nonsmooth convex-concave minimax problem where $\varphi$ and $\psi$ are indicator functions. We select projections of three different closed convex cones to report the performance of the algorithm. In this section, all numerical experiments are implemented by MATLAB R2019a on a laptop with Intel(R) Core(TM) i5-6200U 2.30GHz and 8GB memory.

\subsection{Generalized Absolute Value Equations}\label{sub6.1}
In this part, we focus on generalized absolute value equations \eqref{GAVE} and prove that GAVE can be translated into a smooth linearly constrainted
convex-concave minimax problem. Assume that there exists some $x\in \mathfrak{R}^{n}$ satisfying GAVE \eqref{GAVE}. So far, three techniques have been used to numerically solve GAVE. One is transforming \eqref{GAVE} into a concave minimization to obtain a numerical solution, such as \cite{Mangasarian1997,Mangasarian2014}. The second  technique for solving \eqref{GAVE} is using a generalized Newton method \cite{Mangasarian2008}. The third one is using the matrix splitting technique to solve \eqref{GAVE} (see \cite{Wang2019,Zhou2021}). However, these methods either only obtain the convergence with high probability, or only solve GAVE with special structures, such as $n=m$. Different from these numerical algorithms, Algorithm \ref{alg-30} can be used to solve GAVE by transforming \eqref{GAVE} into a convex-concave minimax problem, and find an $\epsilon$-stationary point without any special structure.

In \cite{Mangasarian2014}, Mangasarian  built the equivalence between GAVE and  the linear complementarity problem, that is, \eqref{GAVE} can be rewritten in the following form
\begin{equation}\label{LCP}
  A(x^+-x^-)+B(x^++x^-)=b,\ \ 0\leq x^+\perp x^-\geq0,
\end{equation}
where $x^+:=(\max\{0,\ x_1\},\ldots,\max\{0,\ x_n\})\in \mathfrak{R}^{n}_+$ and $x^-:=(\max\{0,\ -x_1\},\ldots,\max\{0,\ -x_n\})\in \mathfrak{R}^{n}_+$. Then, we can obtain the solutions of \eqref{LCP} by solving the following linearly constrained minimization
\begin{equation*}
  \begin{array}{ll}
\displaystyle \min_{x^+,\ x^- \in {\Re^n_+}}& \langle x^+,x^-\rangle\\[2pt]
\mbox{subject to} & A(x^+-x^-)+B(x^++x^-)=b,
\end{array}
\end{equation*}
which, by the duality theorem of linear programming, is equivalent to solving the saddle points of the following minimax problem
\begin{equation*}
  \begin{array}{ll}
\displaystyle \min_{x^+\in {\Re^n_+}}\max_{y\in \Re^m}& (b-(A+B)x^+)^Ty\\[2pt]
\mbox{subject to} & (B-A)^Ty\leq x^+.
\end{array}
\end{equation*}
Obviously, the above problem can be expressed as a linearly constrained minimax problem as follows
\begin{equation}\label{LCMMV}
  \begin{array}{ll}
\displaystyle \min_{x^+\in {\Re^n_+}}\max_{z\in {\Re^m_+},y\in \Re^m}& (b-(A+B)x^+)^Ty\\[2pt]
\mbox{subject to} & (B-A)^Ty+z= x^+.
\end{array}
\end{equation}
Then, the proximal gradient multi-step ascent decent method can be used for \eqref{LCMMV} as follows.
\begin{algorithm}\label{alg-V}
Input $(x^0,y^0,z^0,\lambda^0) \in \Re^n \times \Re^n \times \Re^m\times \Re^m$, $\alpha_x >0$, $\alpha_y>0,\ \alpha_z >0$, positive integers $T>1$ and $N>1$ \\
\mbox{}\hspace{0.2cm} {\bf for} $t=0,1,2,\ldots,T-1$, {\bf do}\\[4pt]
\mbox{}\hspace{0.5cm} Set $y^{[0]}(t)=y^t$,\ $z^{[0]}(t)=z^t$\\
\mbox{}\hspace{0.8cm} {\bf for} $k=0,1,2,\ldots, N-1$, {\bf do}\\[4pt]
\mbox{}\hspace{1.2cm} $y^{[k+1]}(t)=y^{[k]}(t)+\alpha_y \left(b - (A+B)^Tx^t+(B-A)\lambda^t\right)$\\[4pt]
\mbox{}\hspace{1.2cm} $z^{[k+1]}(t)=\Pi_{\Re^m_+}(z^{[k]}(t)+\alpha_z \lambda^t)$\\[4pt]
\mbox{}\hspace{0.8cm} {\bf end for}\\[4pt]
\mbox{}\hspace{0.8cm}Set  $x^{t+1}=\Pi_{\Re^n_+}(x^t+\alpha_x \left((A+B)y^{t+1}+\lambda^t\right))$\\[4pt]
\mbox{}\hspace{1.2cm} $\lambda^{t+1}=\lambda^t+\alpha_x [x^{t+1}-(B-A)^Ty^{t+1}-z^{t+1}]$\\[4pt]
\mbox{}\hspace{0.2cm} {\bf end for}\\
\mbox{}\hspace{0.2cm} Return $(x^{t+1},y^{t+1},z^{t+1},\lambda^{t+1})$ for $t=0,1,2,\ldots, T-1$.
\end{algorithm}

Observed that $\psi(y)\equiv0$ and $h(y)=-b^Ty$ in \eqref{LCMMV}, then \textbf{A3} of Assumption  \ref{assum-1e0}  is violated in Theorem \ref{th-main-1e0}, which implies that the iterative complexity of Algorithm \ref{alg-V} is difficult to be proved. However, the performance of Algorithm \ref{alg-V} is the best in numerical experiments. We test three examples of GAVE which are difficult to solve by the previous algorithms. In the first example, $A,\ B$ are symmetric square matrixes and $A+B,\ A-B$ are nonsingular matrixes. The second one is more difficult questions where $A,\ B$ are asymmetric square matrixes and $A+B,\ A-B$ are singular matrixes. Finally, we consider more general linear questions where $A,\ B$ are not square matrixes.
\begin{description}
  \item[(a)] $A,B\in \mathfrak{R}^{n\times n}$ symmetric and $A+B,A-B$ nonsingular
      $$A=\left(
            \begin{array}{ccc}
              1 & 1 & 1 \\
              1 & 0 & 1 \\
              1 & 1 & 1 \\
            \end{array}
          \right),\ B=\left(
            \begin{array}{ccc}
              -1 & 1 & 0 \\
              1 & 2 & 1 \\
              0 & 1 & 1 \\
            \end{array}
          \right),\ b=(-1, 4, 1)^T,
      $$
      Optimal solutions $x^*=(1,-1,-1)$ or $(-1,-1,1)$,\\
      Initial points $\lambda^0=(0,0,0),$\\ $x^0=[0.648679262048621,0.825727149241758,-1.01494364268014], $ \\ $y^0=[-0.471069912683167,0.137024874130050,-0.291863375753573],$\\
      $z^0=[0.301818555261006,0.399930942955802,-0.929961558940129].$
  \item[(b)] $A,B\in \mathfrak{R}^{n\times n}$ asymmetric and $A+B,A-B$ singular
      $$A=\left(
            \begin{array}{ccc}
              -0.5 & 0.5 & 1 \\
              0 & 0.5 & 0.5 \\
              0.5 & 1 & 0 \\
            \end{array}
          \right),\ B=\left(
            \begin{array}{ccc}
              -0.5 & 0.5 & 0 \\
              -1 & 0.5 & 0.5 \\
              0.5 & 1 & 0 \\
            \end{array}
          \right),\ b=(1, 1, 3)^T,
      $$
      Optimal solutions $x^*=(3-2a,a,4-3a)$ with $0\leq a\leq4/3$,\\
      Initial points are the same as in (a).
  \item[(c)] $A,B\in \mathfrak{R}^{m\times n}$ with $m\neq n$
      $$A=\left(
            \begin{array}{cccc}
              0.01  & \cdots & 0    &  \\
               &   \ddots &  &  { \Large\textbf{0}_{100\times100}}   \\
               &   0.01 & 0   &   \\
               &   0.01 & 0   &  \\
            \end{array}
          \right)^T\in \mathfrak{R}^{200\times 100}
      ,$$
      $$\ B=\left(
            \begin{array}{rrrrc}
              0.01  &-1& \cdots & -1    &  \\
               &   \ddots & & &  { \Large-\textbf{1}_{100\times100}}   \\
               &  & -0.99 & -1   &   \\
               &  & 0.01 & -0.99   &  \\
            \end{array}
          \right)^T\in \mathfrak{R}^{200\times 100},$$
          $$b=(-7.99,-7.01,-6, \cdots,-6,-6.01, -5)^T\in\mathfrak{R}^{200},
      $$
      Initial points $x^0=y^0=z^0=\lambda^0=(0,\cdots,0)\in\mathfrak{R}^{200}.$
\end{description}

We denote $\textbf{0}_{100\times100}$ and $\textbf{1}_{100\times100}$ as matrices with all components of 0 and 1, respectively. In Table \ref{tab1}, we report the CPU time and iteration  obtained by running Algorithm \ref{alg-V} for solving the above three problems. Even for the general matrix form of GAVE (c), Algorithm \ref{alg-V} can quickly converge to the solution of GAVE. Moreover, for GAVE with infinite solutions (b), Algorithm \ref{alg-V} can also find the approximate solution of GAVE within a small error.


\begin{table}[http]
\begin{center}
\begin{tabular}{|c|c|c|c|c|c|}
\hline
Problem & $\alpha_x=\alpha_y=\alpha_z$ & N                            & T &  Time (s)& Error \\ \hline
(a)     & 0.05                & 5       & 119&  1.66&  8.66e-5       \\ \hline
(b)     & 0.01                & 40                        & 46  & 0.19 & 2.34e-2      \\ \hline
(c)     & 0.01                & 5                        & 3  & 7.42e-2  & 3.18e-12      \\ \hline
\end{tabular}
\caption{ Numerical results of Algorithm \ref{alg-30} for GAVE \eqref{GAVE}. ($\alpha_x$- step sizes w.r.t $x$; $\alpha_y$- step sizes w.r.t $y$; $\alpha_z$- step sizes w.r.t $z$; $N$- number of the inner loop; T-number of the outer loop; Time-CPU time; Error-$\|Ax+B|x|-b\|$.)}\label{tab1}
\end{center}
\end{table}

\subsection{Linear Regression }\label{sub6.2}
In this part, we consider the well-known linear regression problem with joint linearly constraints as follows
\begin{equation}\label{lr}
\begin{array}{rl}
\displaystyle \min_{x \in \Re^n}\max_{y \in \Re^m} & f(x,y)=\frac{1}{m}\left[-\frac{1}{2}\|y\|^2-b^Ty+y^TKx\right]+\frac{\lambda}{2}\|x\|^2\\[10pt]
\mbox{subject to} & Ax+By+c=0_p,
\end{array}
\end{equation}
where the rows of the matrix $K\in \Re^{m\times n}$,\ $A\in \Re^{p\times n}$ and $B\in \Re^{p\times m}$ are generated by a Gaussian distribution $\mathcal {N}(0,I)$. In the following experiments, let $n=m$, $b=0$, $c=0$ and $\lambda=1/m$.

For solving \eqref{lr},
in Algorithm \ref{alg-30}, we set the step sizes $\alpha_x=0.3,\ \alpha_y=1$ and the number of inner loops is selected as $N=3$ and randomly select the initial points. Algorithm \ref{alg-30} is terminated when at $t$-iteration, $(x^{t},y^{t},\lambda^{t})$ is an  $\epsilon$-KKT point; {\it i.e.}, at $t$-iteration, $
G^{f,\varphi}_L(x^{t},y^{t},\lambda^{t})\leq10^{-7},
G^{f,0}_L(x^{t},y^{t},\lambda^{t})\leq10^{-7},
G^{f,\psi}_L(x^{t},y^{t},\lambda^{t})\leq10^{-7}.
$

\begin{figure}[htbp]
  \centering
  \subfigure[$n=10$]{
    \label{fig:subfig:10} 
    \includegraphics[width=2.9in]{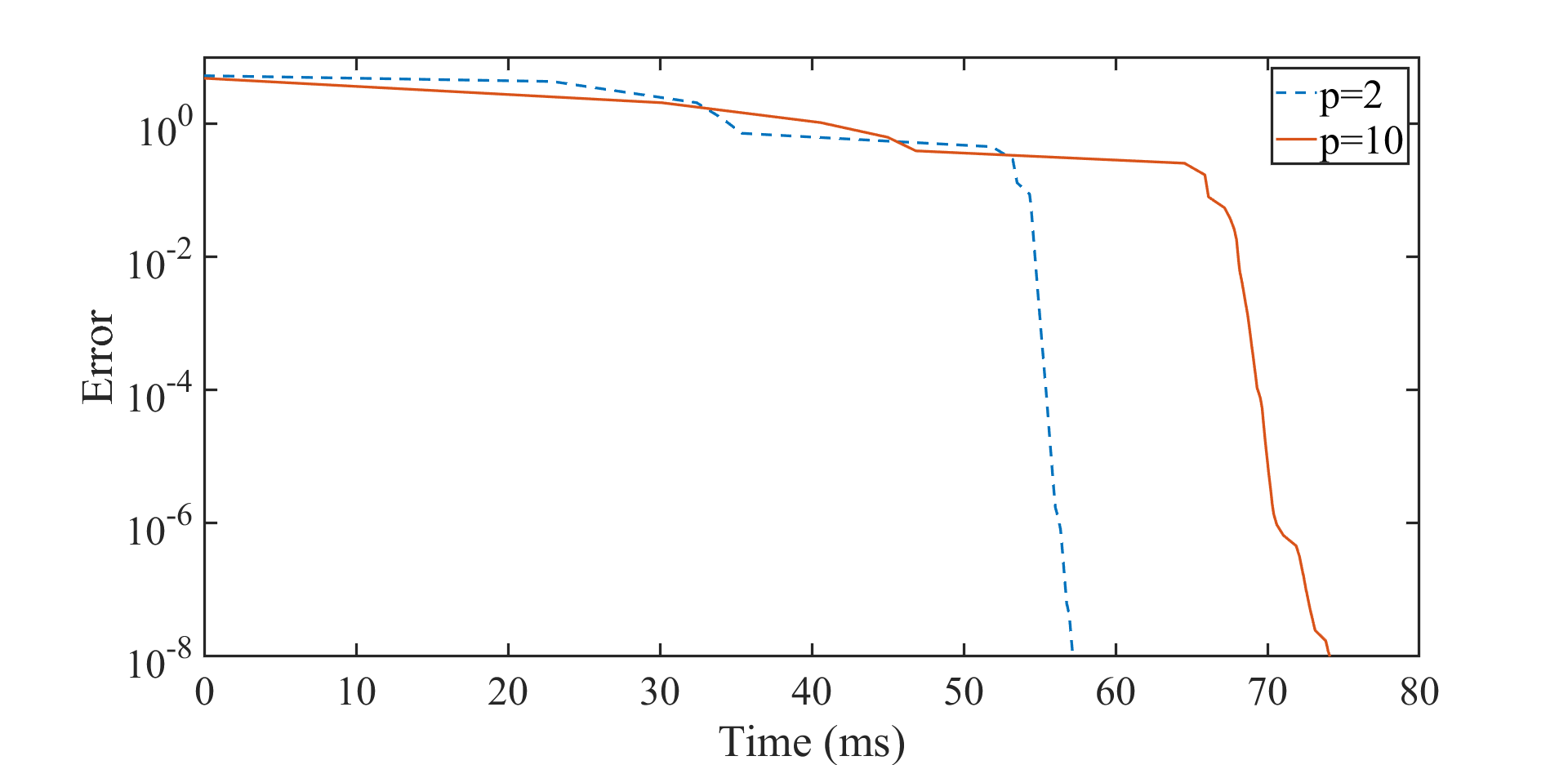}}
  \subfigure[$n=100$]{
    \label{fig:subfig:100} 
    \includegraphics[width=2.9in]{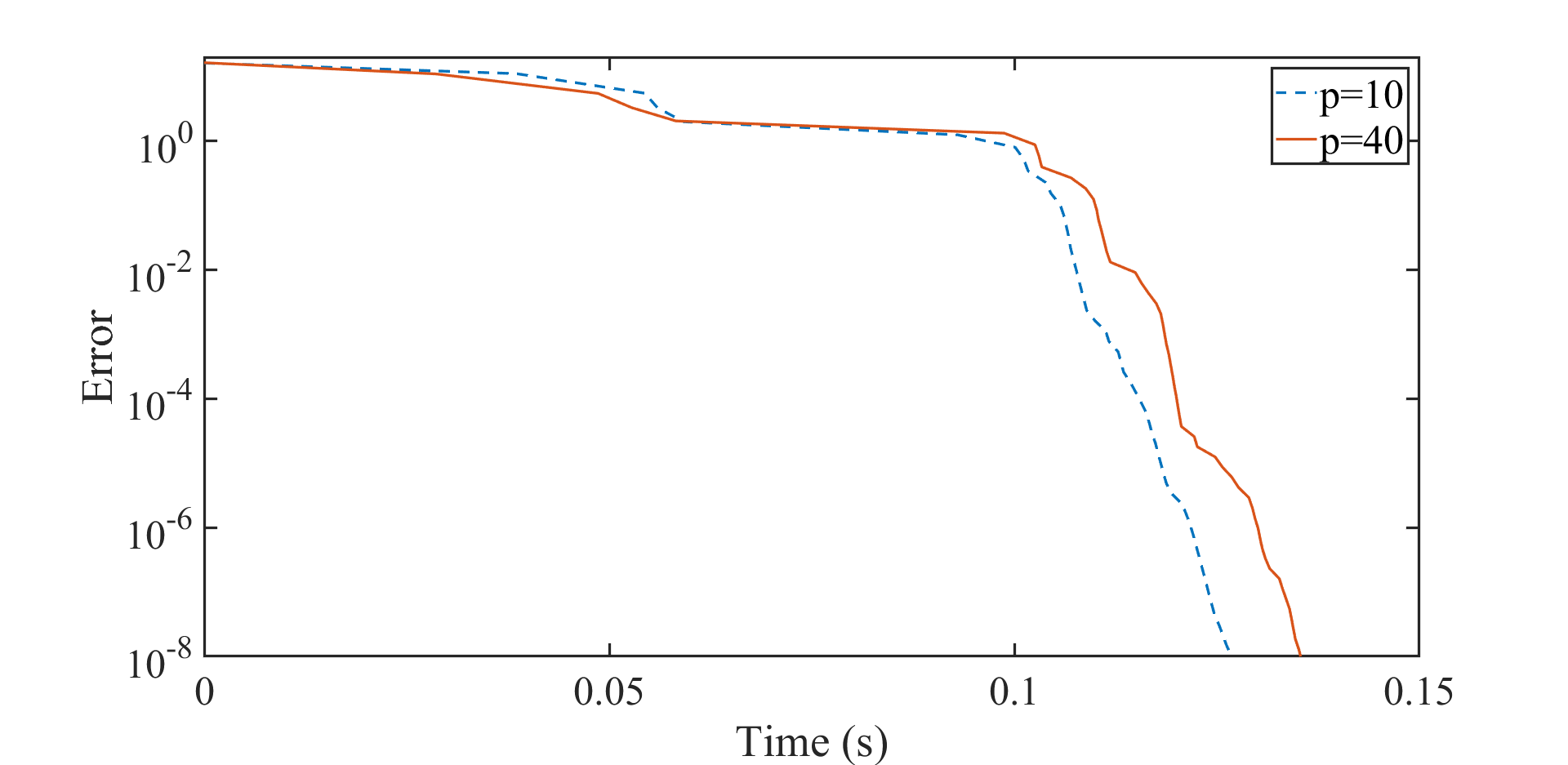}}

    \subfigure[$n=1000$]{
    \label{fig:subfig:1000} 
    \includegraphics[width=2.9in]{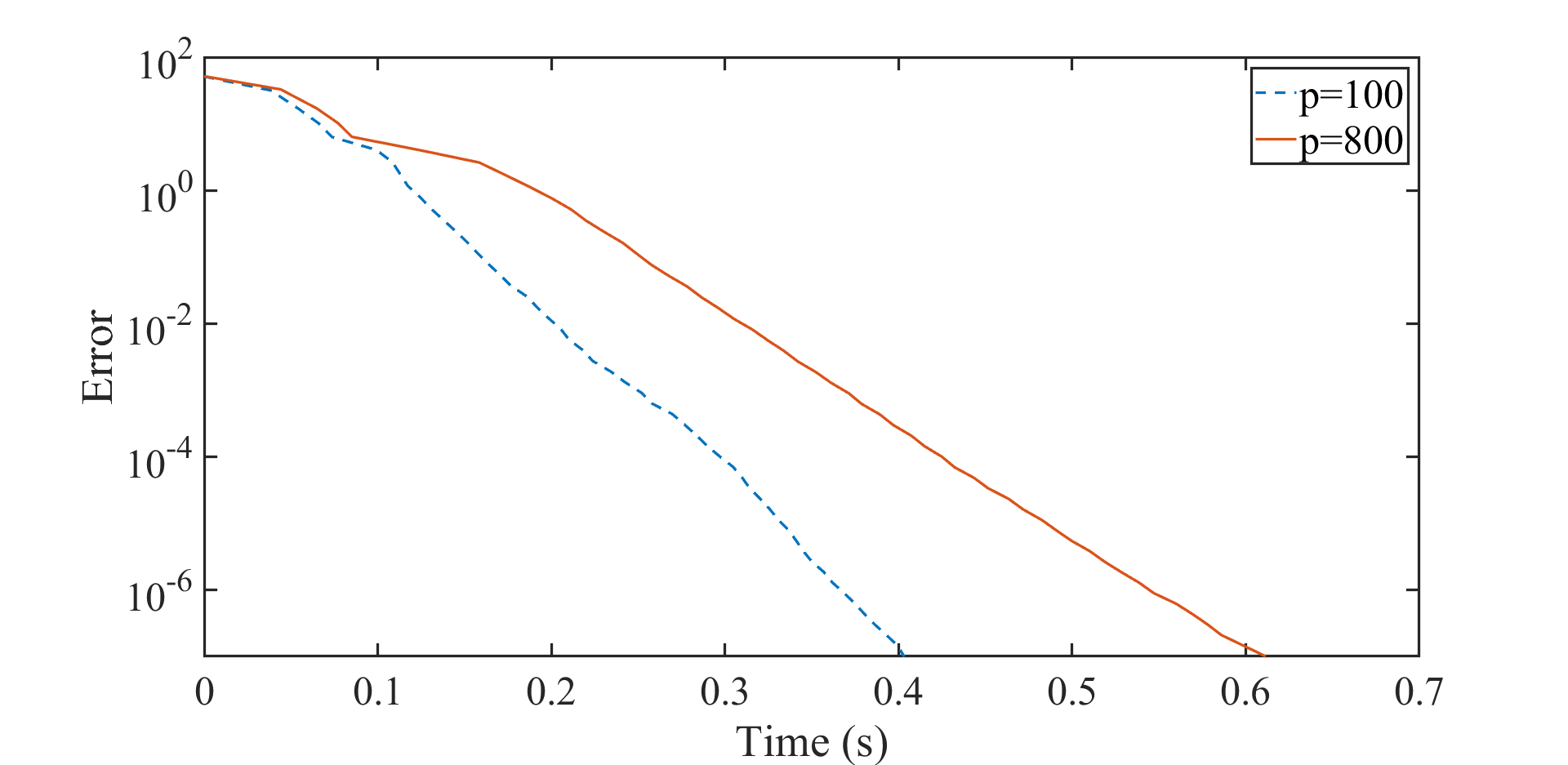}}
    \subfigure[$n=10000$]{
    \label{fig:subfig:1000} 
    \includegraphics[width=2.9in]{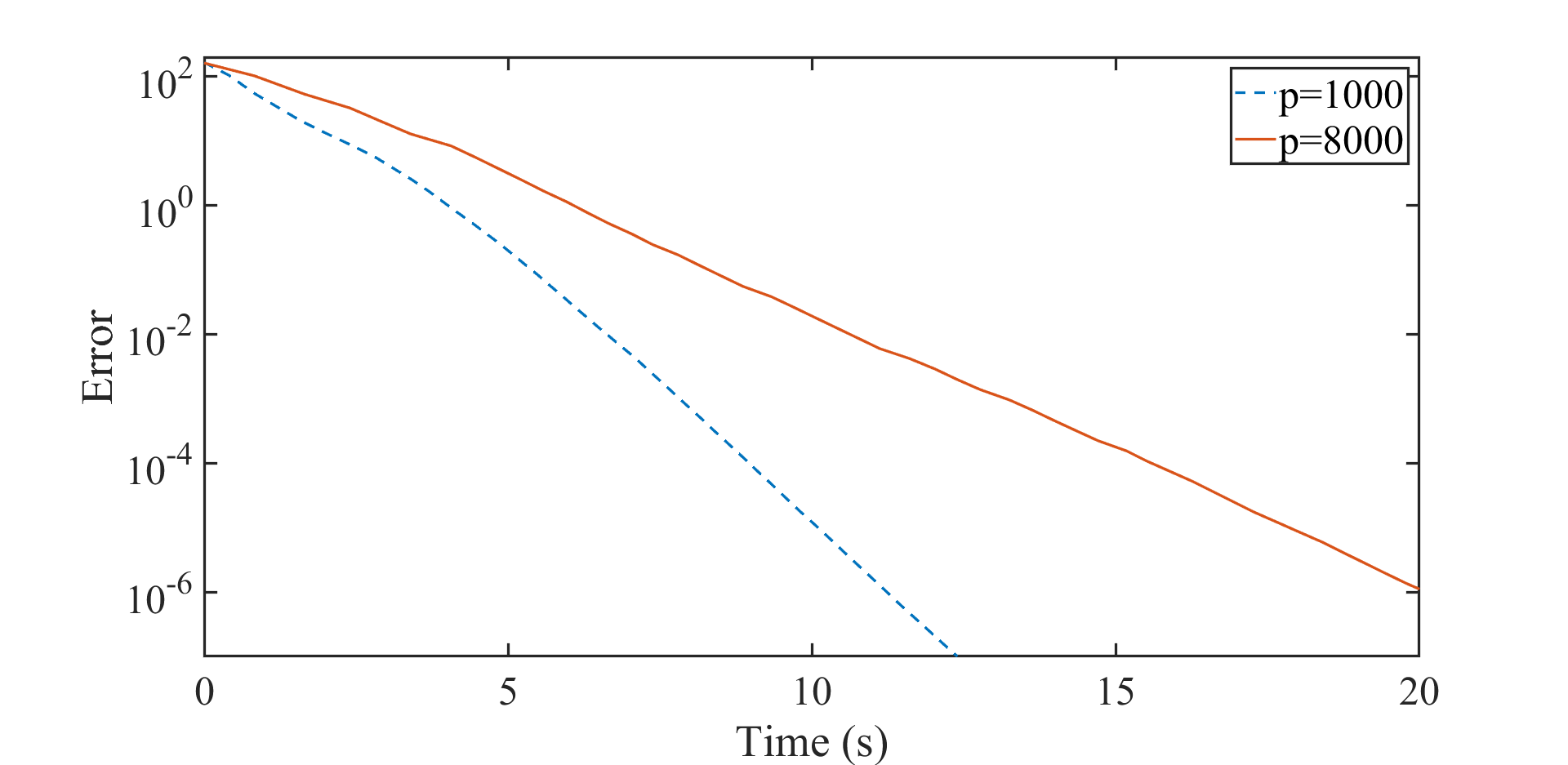}}
  \caption{The trend of the error for linear regression problems \eqref{lr} with respect to CPU time. (Error-$
\|G^{f,\varphi}_L(x^{t},y^{t},\lambda^{t})+
G^{f,0}_L(x^{t},y^{t},\lambda^{t})+
G^{f,\psi}_L(x^{t},y^{t},\lambda^{t})\|.
$)}
  \label{fig:subfig3} 
\end{figure}

Figure \ref{fig:subfig3} shows  the performance of Algorithm \ref{alg-30} for \eqref{lr} with different dimensions and constraints.  As seen
from Figure \ref{fig:subfig3}, Algorithm \ref{alg-30} converges linearly and ensures that the optimal solution can be found quickly (within 30 seconds) for large dimension ($n=O (10^4)$).

\subsection{Generalized Linear Projection Equations}\label{sub6.3}
In this subsection, we consider the generalized linear projection
equations of the form
\eqref{glpe}. One of classic problems of \eqref{glpe} is linear projection
equations (e.g. \cite{Eshaghnezhad2016,Heemels2001,XiaWang2000}) of the form:
\begin{equation}\label{lqe}
P_K(u-Mu-q)=u,
\end{equation}
where $M\in \mathfrak{R}^{n\times n},\ q\in \mathfrak{R}^{n}$ and $K$ is a closed polyhedron. \eqref{lqe} can be expressed as the projection equations in the following standard form:
$$(I-M)^{-1}x-x_K=-q,\ {\rm s.t.}\ x_K=P_K(x).$$
There are some well-established techniques used to solve \eqref{lqe}, such as Lemke's methods \cite{Miller1991},  interior-point methods \cite{Castronuovo2001} and neural networks \cite{XiaWang2000}. However, there are few numerical algorithms for developing generalized linear projection equations \eqref{glpe}. In this part, we propose a new technique for solving generalized linear projection equations, where \eqref{glpe} is converted to a minimax problem \eqref{minimax} and solved with Algorithm \ref{alg-30}.

In the following, assume that $K$ is a closed convex cone. Hence, for any $x\in \mathfrak{R}^{n}$, it can be expressed as
$$x=P_K(x)+P_{K^{\circ}}(x),$$
where $K^{\circ}$ represents the polar cone of $K$. Then, \eqref{glpe} can be rewritten as
\begin{equation*}
  \begin{array}{ll}
\displaystyle \min_{x_{K}\in {K},x_{K^{\circ}}\in {K^{\circ}}}& \langle x_{K},x_{K^{\circ}}\rangle\\[2pt]
\mbox{subject to} & A(x_{K}+x_{K^{\circ}})+Bx_{K}=b,
\end{array}
\end{equation*}
where $x_{K}=P_K(x), x_{K^{\circ}}=P_{K^{\circ}}(x)$.
By the duality theorem of linear programming, the above problem is equivalent to the following linearly constrained minimax problem
\begin{equation}\label{dglpe}
  \begin{array}{ll}
\displaystyle \min_{x_{K}\in {\Re^n}}\max_{y\in \Re^m,z\in {\Re^n}}& \delta_K(x_{K})+(b-(A+B)x_{K})^Ty-\delta_{K^{\circ}}(z)\\[2pt]
\mbox{subject to} & A^Ty+z= x_{K},
\end{array}
\end{equation}
where $y$ is a Lagrange multiplier and $z$ is an auxiliary variable, and $\delta_K(x)$ is an indicator function with $\delta_K(x)=0$ for $x\in K$ while $\delta_K(x)=+\infty$ for $x\notin K$. The solutions to equations \eqref{glpe} can be found by applying Algorithm \ref{alg-30} to solve \eqref{dglpe}.

We tested projection equations \eqref{glpe} with symmetric and asymmetric cones, which are the positive octant cone, second-order cone and 1-norm cone, respectively expressed as
\begin{itemize}
  \item Positive octant cone (symmetric): $K=\Re^n_+$;
  \item Second-order cone (symmetric): $K=\{x:=(s_0;\bar{s})\in\Re^n:\|\bar{s}\|\leq s_0,\bar{s}\in\Re^{n-1}\}$, where $\|\cdot\|$ is the Euclidean norm;
  \item 1-norm cone (asymmetric): $K=\{x:=(s_0;\bar{s})\in\Re^n:\|\bar{s}\|_1\leq s_0,\bar{s}\in\Re^{n-1}\}$, where $\|\cdot\|_1$ is the 1-norm.
\end{itemize}
Let nonsingular asymmetric matrices $A,B\in \mathfrak{R}^{n\times n}$ be given as follows.
      $$A=\left(
            \begin{array}{ccccc}
              -1 & 0 & 1 & 0 & 0 \\
              1 & 0 & -1 & 1 & 1\\
              -1 & 1 & 1 & 0 & 0 \\
              0 & 1 & 1 & -1 & 0 \\
              1 & -1 & 1 & 0 & 1
            \end{array}
          \right),\ B=\left(
            \begin{array}{ccccc}
              0.5 & 0.5 & 1 & 0 & -1 \\
              1 & 0 & 0.5 & 1 & 2\\
              1 & -1 & 1 & 0.5 & 1 \\
              0 & 0 & -1 & -0.5 & 1 \\
              1 & 0 & 0 & 0 & 0.5
            \end{array}
          \right),\ b=(6.5,    5,    8.5,   -1.5,    8.5)^T.
      $$
Algorithm \ref{alg-30} is used to solve \eqref{dglpe} of the above examples. In numerical experiments, we set the step sizes $\alpha_x=\alpha_y=\alpha_z=1/|{\rm det}(A+B)|$, and the number of inner loops is selected as $N=5$, and randomly select the initial points. Algorithm \ref{alg-30} is terminated when at $t$-iteration, $\|Ax^t+Bx^t_K-b\|\leq10^{-14}$.

\begin{figure}
  \centering
\includegraphics[width=5in]{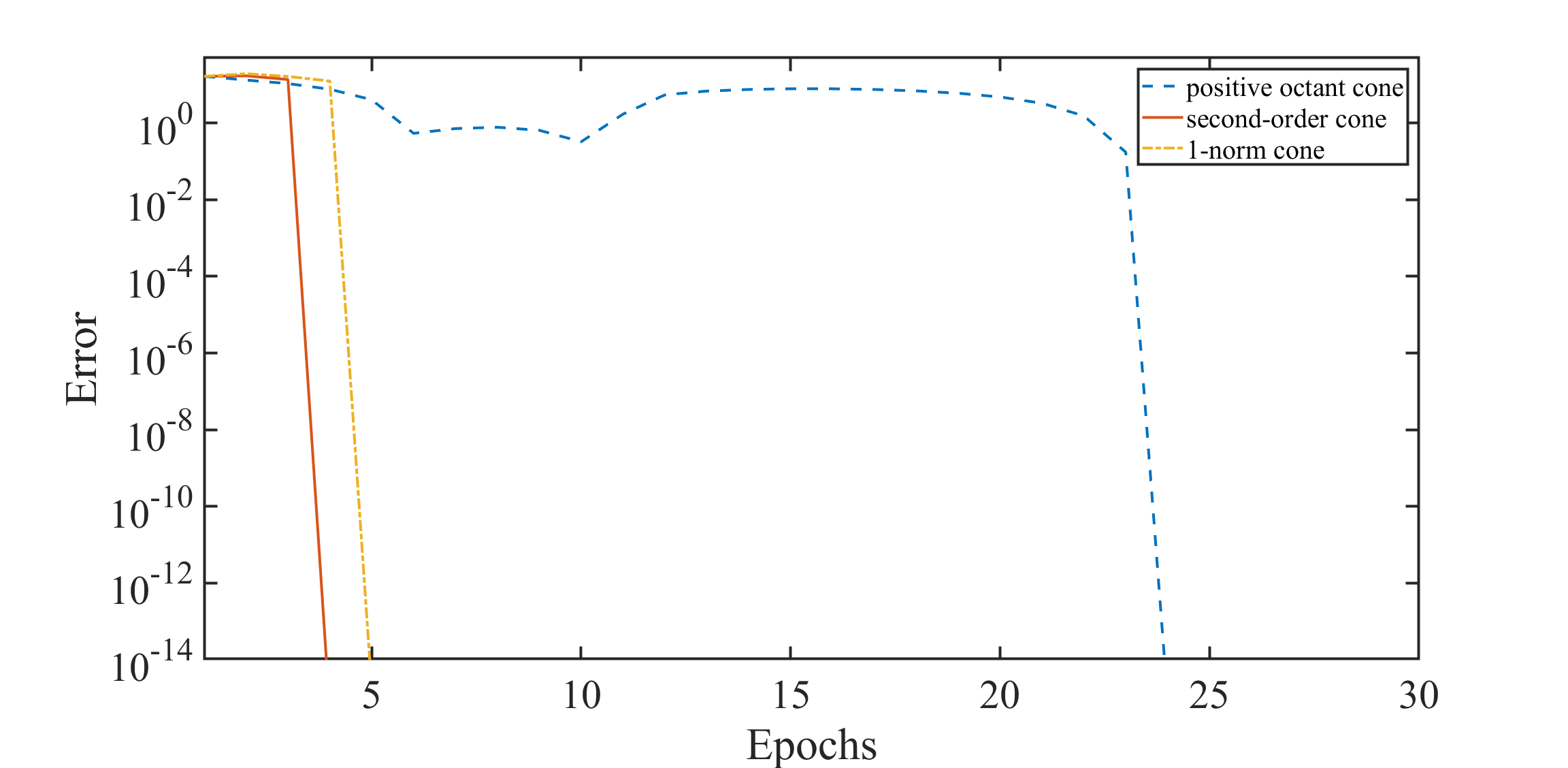}
  \caption{The trend of the error for generalized linear projection equations \eqref{glpe} under different cones. (Error-$\|Ax+Bx_K-b\|$)}
  \label{fig:subfig4} 
\end{figure}

The effect of Algorithm \ref{alg-30} for \eqref{glpe} is given in Figure \ref{fig:subfig4}. For all symmetric and asymmetric cone, Algorithm \ref{alg-30} converges rapidly and obtains nearly accurate solutions of \eqref{glpe}.
\section{Some Concluding Remarks}\label{Sec4}
 \setcounter{equation}{0}
\quad \, Nonsmooth linearly constrained minimax optimization problems are an important class of optimization problems with many applications, such as  generalized absolute value equations and linear regression problems. However, there are few  numerical algorithms for solving this type of problems if there are joint linear constraints. We developed a conceptual alternating coordinate ascent-decent method in which the global convergence is guaranteed if the subproblems with respect to $y$ are solved sufficiently accurately. Combining specific numerical methods (proximal gradient methods) to solve the subproblems, we proposed a proximal gradient multi-step ascent decent method and demonstrated the iteration complexity bound for an $\epsilon$-stationary point in ${\cal O}\left(\epsilon^{-2}\log  \epsilon^{-1}  \right)$ iterations under mild conditions. Finally, we applied the proximal gradient multi-step ascent decent method to generalized absolute value equations and linear regression problems and proved the linear convergence of the method.

There are many interesting problems worth considering. Currently, the iteration complexities for  the proximal gradient multi-step ascent decent method  is obtained under Assumption \ref{assum-domain-B} or Assumption \ref{assum-3e0}. How to guarantee Assumption \ref{assum-3e0}  is an important issue for further study.

\end{document}